\documentclass[a4paper, fleqn]{article}
\usepackage[english]{babel}
\usepackage{graphicx}
\usepackage{amsmath, amssymb, graphics, float, fullpage}
\usepackage{amsthm}
\usepackage{makeidx}
\usepackage{lipsum}
\usepackage{tikz}

\numberwithin{equation}{section}

\newtheorem{theorem}{Theorem}[section]
\newtheorem{lemma}[theorem]{Lemma}
\newtheorem{proposition}[theorem]{Proposition}
\newtheorem{corollary}[theorem]{Corollary}
\newtheorem{definition}[theorem]{Definition}
\newtheorem{remark}[theorem]{Remark}

\newtheorem{assumption}[theorem]{Assumption}

\providecommand\phantomsection{}

\makeindex

\begin{document}

\clearpage
\thispagestyle{plain}

\title{THE MATCHING CONDITION FOR LARGER SIZE RIEMANN-HILBERT PROBLEMS}
\author{L. D. Molag}
\maketitle
\begin{center}
KU Leuven, Department of Mathematics,\\
Celestijnenlaan 200B box 2400, BE-3001 Leuven, Belgium.\\
E-mail: leslie.molag@kuleuven.be
\vspace{1.0cm}

\end{center}

\begin{abstract}
In a larger size Riemann-Hilbert problem matching the local parametrices with the global parametrix is often a major issue. In this article we present a result that should tackle this problem in natural situations. We prove that, in a general setting, it is possible to obtain a double matching, that is, a matching condition on two circles instead of one circle. We discuss how this matching approach can be used to obtain local scaling limits of correlation kernels and apply our result to several examples from the existing literature. 
\end{abstract}

\tableofcontents

\section{Introduction and statement of results}

\subsection{Introduction}

A Riemann-Hilbert problem (RHP) asks for a complex function, possibly matrix valued, that is analytic outside an oriented contour, has certain jump properties on this contour and has a specified asymptotic behavior at infinity and at other predescribed points. 

Nowadays $n$-dependent RHPs are used as an effective tool to find the large $n$ behavior of functions of interest. This is mainly a consequence of the influential papers \cite{FoItKi, DeZh}. In these papers a formalism, commonly refered to as a Riemann-Hilbert analysis, is developed in which one can keep transforming a RHP to a new one, eventually reaching a RHP that is trivial enough to draw conclusions about the large $n$ behavior (e.g., see \cite{KaMcMi} or \cite{BlLi}). This method is known as the \textit{Deift-Zhou steepest descent analysis}. A particularly nice insight is the correspondence between RHPs and orthogonal polynomials \cite{De}. 

Though RHPs have been very effective in the $2\times 2$ case, larger size RHPs often have one problematic step. This step is the so-called \textit{matching}. After transforming a RHP several times, one usually finds oneself in a situation where the RHP at hand has to be approximated by a so-called global parametrix, which is defined away from certain points, and local parametrices, which are defined in a small neighborhood of these points. We shall refer to these points as the \textit{special points} henceforth. To obtain the matching condition one should modify the local parametrices in such a way that their values on the boundary of their domain are very close to the values of the global parametrix there. As it turns out, it is generally hard to find this modification when we have an $m\times m$ matrix valued RHP with $m>2$ (see Section \ref{sec:difficulty}). 

In the literature one finds several methods to obtain this matching for specific larger size RHPs \cite{BeBo, DeKu, DeKuZh, KuMFWi}, but no general result was known. In particular, the somewhat ad hoc and technical nature of these methods seems to become a problem when $m$ is large. What is also important, is that we sometimes want to consider RHPs of arbitrary size rather than a fixed size. In \cite{KuLo} an analysis for a RHP of arbitrary size was carried out, but there the matching effectively boiled down to a decomposition in $2\times 2$ matching problems. In a recent article \cite{KuMo} by Kuijlaars-Molag an iterative method was developed for a specific $3\times 3$ RHP to obtain the matching, which they expected to be applicable to other larger size RHPs as well. Indeed this is the case. Very recently, the method was succesfully adopted in \cite{SiZh}. However, in this paper we develop the method further considerably. We prove in Section \ref{sec:proofMainResult} that, assuming a natural situation, one can modify the global and local parametrix in a local way such that one has a \textit{double matching}, that is, a matching on two circles instead of one. This main result is formulated in Theorem \ref{lem:matching}

In many applications a RHP aims to elucidate the large $n$ behavior of functions of interest, e.g., orthogonal polynomials, oscillatory integrals and correlation kernels. Such behavior is frequently expressed  in the form of a scaling limit. We clarify in Section \ref{section:scalingLimits} how the double matching can be used to obtain scaling limits when the RHP aims to find the large $n$ behavior of an associated correlation kernel. This is contained in Theorem \ref{thm:corKerScaling}. It is to be noted that Lemma \ref{thm:scalingb} is not limited to the case of correlation kernels, and it may be useful for other type of scaling limits as well. 

We conclude with applications of our result to several examples from the existing literature, which is the topic of Section \ref{sec:examples}. All these examples concern scaling limits of correlation kernels. The first example is worked out in detail, and should provide the reader with some general intuition of how to apply the main result. None of these examples lead to new results, but the fact that the main result can be applied to all of them gives credit to its applicability, and strengthens my belief that it will be a useful tool for future RH analyses. 

I would like to point out that, with the addition of some convergence properties, the conclusion of the main result seems to hold even when $m=\infty$. Such infinite size RHPs, and even more general operator valued ones, have been considered (e.g., see \cite{ItKo2}), but for the moment matching global with local parametrices does not appear to be an issue. 

As a final remark, let me say that there are also examples where a construction of the local parametrix with special functions did not seem to be possible and instead a small norms argument had to be used (e.g., see \cite{KrMc}). It is possible then, that one is not in a sufficiently generic situation, and perhaps \text{Theorem \ref{lem:matching}} cannot be applied without doing some extra work. On the other hand, if one would have to shrink the boundary circle in order to apply a small norms argument, then Theorem \ref{lem:matching} might actually be useful under the right circumstances even for non-standard $2\times 2$ RHPs.

\pagebreak

\subsection{Set up} \label{sec:setUp}

Here and in the rest of the paper we use the notation $D(0,r)$ for the disc $|z|<r$ and we orient its boundary positively. We also use the notation $A(0;r_1,r_2)$ for the annulus $r_1<|z|<r_2$. Throughout the article $m$ is a fixed positive integer. By a contour $\Sigma$ we shall mean a finite union of smooth curves, that intersect in at most a finite number of points. $\Sigma^0$ will denote $\Sigma$ minus the intersection points. 

Given an oriented contour $\Sigma$ in $\mathbb C$ and a possibly $n$-dependent function $v:\Sigma\to\mathbb C^{m\times m}$ (the jump matrix), the associated RHP asks for a function $Y:\mathbb C\setminus\Sigma\to\mathbb C^{m\times m}$ such that
\begin{itemize}
\item[RH-Y1] $Y$ is analytic.
\item[RH-Y2] $Y$ has boundary values $Y_\pm(z)$ for $z\in \Sigma^0$ that satisfy the jump relation
\begin{align} \nonumber
Y_+(z) = Y_-(z) v(z).
\end{align}
Here the $+$ or $-$ determines from which side we approach the contour. 
\item[RH-Y3] $Y$ has a specified asymptotic behavior near $\infty$ and near other predescribed points.
\end{itemize}
As is common practise, we do not indicate the $n$-dependence of $Y$, or any of the transformations thereof. We are being intentionally vague in the formulation of RH-Y3, as it shall not be of any relevance in this paper, though we mention that the asymptotic behavior in RH-Y3 may also depend on $n$. The usual approach is to transform $Y$ according to
\begin{align} \nonumber
Y \mapsto X \mapsto T \mapsto S,
\end{align}
where we have consecutive transformations that we call the \textit{first transformation}, the \textit{normalization} and \textit{opening of lenses}.  Sometimes one deviates slightly from the particular transformations used, but one basically always ends up with a RHP for some function $S$ that has to be approximated by appropriately chosen functions. Namely, for each special point it should be approximated by a function $P$ on some disc around the special point. In the remaining part, i.e., the complex plane minus the discs around all the special points, $S$ should be approximated by a function $N$. The first approximating function is called the local parametrix (at the particular special point), it is usually constructed using specific special functions. The second is called the global parametrix, it generally has the interpretation of a limiting behavior of $S$ as $n\to\infty$ away from the special points and as such it should be an $n$-independent function, although $n$-dependent global parametrices can arise in multi-cut cases (e.g., see \cite{DeItZh}). In general, we have an expression for $N$ for the full complex plane minus the special points (although a choice has to be made for the values on the jump contour). It is generally clear from the RH analysis what points are special and what points are not. 

In the remainder of this paper we assume that there is a special point at the origin. There may be other special points as well, but we will focus solely on $z=0$. Notice that there is no loss in generality in doing so, since we can always translate $z$ without non-trivially changing the RHP for $S$. It is also reasonable to assume that the set of special points is discrete. We thus include the following assumption. 

\begin{assumption} 
The local parametrix $P$ will be defined on a closed disc $\overline{D(0,r)}$ for some $r>0$ and there are no other special points in the disc. $N$ will be a function on $\mathbb C$ minus the special points.
\end{assumption}

In general $P$ and $N$ are analytic, except on the jump contour (together with the end points). In a natural situation, $P$ has the same jumps as $N$ has, but also some additional jumps that correspond to jumps of $S$ on the lips of lenses.   The \textit{matching condition} demands that there exists a $\delta>0$ such that
\begin{align} \label{eq:matchingCondition}
P(z) N(z)^{-1} &= \mathbb I + \mathcal O(n^{-\delta}), & \text{ uniformly for }|z|=r, \quad\text{ as }n\to\infty. 
\end{align}
It is also allowed to obtain the matching on a shrinking circle (e.g., see \cite{BlKu}), that is, $r$ does not have to be fixed but it may depend on $n$ such that $r\to 0$ as $n\to\infty$. It often happens that our initial construction for $P$, let's call it $\mathring P$, almost works, but that we need a non-singular analytic function $E_n$, called the \textit{analytic prefactor}, such that
\begin{align} \nonumber
E_n(z) \mathring P(z) N(z)^{-1} &= \mathbb I + \mathcal O(n^{-\delta}), & \text{ uniformly for }|z|=r, \quad\text{ as }n\to\infty. 
\end{align}
We would then thus identify $P = E_n \mathring P$ in order to obtain the matching \eqref{eq:matchingCondition}.

When the matching is successful, one applies what is known as the \textit{final transformation}. For convenience, since we focus on the special point $0$, we absorb any possible local parametrices around other special points into the definition of $N$.  One then transforms $S\mapsto R$, according to
\begin{align} \label{eq:finalOriginal}
R(z) = \left\{\begin{array}{ll}
S(z) N(z)^{-1}, & |z|> r,\\
S(z) P(z)^{-1}, & |z|<r.
\end{array}\right.
\end{align}

As one can check, $R$ then has a jump matrix on $\partial D(0,r)$ that behaves uniformly as $\mathbb I + \mathcal O(n^{-\delta})$ as $n\to\infty$. If any possible other jumps are also close to $\mathbb I$ then, under some additional conditions, one may draw strong conclusions about the large $n$ behavior of $R$ via a general theorem (e.g., see \cite{De} or \cite{BlLi}). One may then invert all the transformations of the RHP to obtain the large $n$ behavior of $Y$. The fact that the jump matrices of $R$ need to be close to $\mathbb I$ in order to apply this general theorem, is the reason why one desires to obtain the matching condition in the first place.

\subsection{The issue with the matching for larger size RHPs} \label{sec:difficulty}

To understand why the matching is hard for larger size RHPs, we have to be precise about the structure of the local parametrix. We are assuming here, that all the necessary ingredients are explicit, as opposed to the situation where we have to use a small norms argument, for example. We point out that a detailed example can be found in Section \ref{sec:examples}. It may be beneficial to a reader that is unexperienced in solving larger size RHPs, or a reader that would like to have some explicit handhold, to have a look at that example before proceding with this subsection. 

Let us sketch the general picture when attempting to match the global with the local parametrix. Usually, our initial construction for a solution to the local parametrix problem is of the form
\begin{align} \label{eq:defmathringP}
\mathring P(z) = \mathring E_n(z) \Psi\left(n^b f(z)\right) D(z) e^{n D_\varphi(z)},
\end{align}
where the expressions are as follows. $D_\varphi(z)$ is a diagonal matrix whose components are a linear combination of the so-called $\varphi$-functions that originated from the normalization and opening of lenses. Sometimes one uses so-called $\lambda$-functions instead of the $\varphi$-functions (e.g., see \cite{KuMFWi}). $D(z)$ consists of Sz\H{e}go functions that are used to make the jumps constant, it is often some constant diagonal matrix power of $z$ depending on how trivial the $z$ dependence of the jumps is. The combination $D(z) e^{n D_\varphi(z)}$ has the effect of reducing the initial local parametrix problem to one that has constant jumps. To solve such a problem one considers the same problem but with the curves of the jump contour extended to infinity, dividing $\mathbb C$ into different sectors, which is often called the \textit{model Riemann-Hilbert problem}. Since this terminology is also frequently used to describe the global parametrix problem, I prefer to use \textit{bare parametrix problem} instead, as in \cite{BeBo} (although a trivial dependence on $z$ was allowed there). This bare parametrix problem is then solved by $\Psi$, which is constructed using specific special functions. Well-known examples of such functions are Airy functions \cite{DeKrMcZh}, Bessel functions \cite{KuMcVAVa} or parabolic cylinder functions \cite{DeZh}, but more exotic ones such as for example functions corresponding to Painlev\'e equations \cite{ClKuVa}, and Meijer G-functions \cite{BeBo}, have also been used. In general, $\Psi$ has an asymptotic series behavior
\begin{align} \label{eq:behavPsi}
\Psi(\zeta) \sim \left(\mathbb I + \frac{C_1}{\zeta} + \frac{C_2}{\zeta^2} + \ldots\right) B(\zeta) e^{\theta(\zeta)},
\end{align}
as $\zeta\to\infty$. The coefficients $C_1,C_2,\ldots$ in the asymptotic series are $m\times m$ constant matrices. $\theta(\zeta)$ is some diagonal matrix whose entries (in most cases) are multiples of $\zeta^\frac{1}{b}$. These multiples may be different in different sectors, in the sense that they may be permuted. $B(\zeta)$ is some matrix-valued function with power law behavior. Frequently, $B$ is a constant diagonal matrix power of $\zeta$. Then, viewing one particular sector, the expression \eqref{eq:behavPsi} is known as a formal Birkhoff invariant \cite{Bi}. It is well-known that such an invariant is a formal fundamental solution of a formal meromorphic differential equation, which should correspond to the particular special functions in the construction. $\Psi$ is evaluated in $n^b f(z)$, where $f$ is a conformal map around $0$ that maps $0$ to itself, and, in standard situations, maps positive numbers to positive numbers. $b$ is a positive number that plays a role in scaling limits concerning the solution of our RHP. Lastly, $\mathring E_n$ might be our first guess for an analytic prefactor, which we may also set to $\mathbb I$ if there does not appear to be an obvious candidate. We repeat that explicit examples can be found in \text{Section \ref{sec:examples}}.

On some circle around the origin of fixed radius we then have as $n\to\infty$
\begin{align} \label{eq:matchingAsympS}
\mathring P(z) N(z)^{-1} = \left(\mathbb I+\frac{C_1}{n^b f(z)}+\frac{C_2}{(n^b f(z))^2}+ \ldots + \frac{C_k}{(n^b f(z))^k}+\mathcal O(n^{-c})\right) E(z)^{-1},
\end{align}
for some $c>0$, a non-negative integer $k$ and some $n$-dependent analytic function $E(z)^{-1}$ (the inverse is included for reasons that will become clear shortly). In fact, \eqref{eq:matchingAsympS} is valid on shrinking circles as well, provided that their radii shrink sufficiently slower than order $n^{-b}$ as $n\to\infty$. If one has done the work to write $\Psi$ in the asymptotic form \eqref{eq:behavPsi} then, setting $\mathring E = \mathbb I$, we may take
\begin{align} \label{eq:defmathE}
E(z) = N(z) D(z)^{-1} e^{-n D_\varphi(z)-\theta(n^b f(z))} B(n^b f(z))^{-1}.
\end{align}
One would have to check on a case-by-case basis that $E$ is indeed non-singular and  analytic around $z=0$. In some cases though, it might be easier to arrive at a situation similar to \eqref{eq:matchingAsympS} (with possibly different coefficients $C_1, C_2, \ldots$) by choosing a convenient prefactor $\mathring E$. We will see such cases in Section \ref{sec:examples}.

Heuristically, one is in the situation \eqref{eq:matchingAsympS} because away from the origin the effects of the additional jumps of $\mathring P$, usually corresponding to the lenses, become so weak for large $n$ that $\mathring P N^{-1}$ is approximately an analytic function. The effect of the extra jumps is essentially contained in the $\mathcal O(n^{-c})$ term.

One may take as many terms of the asymptotic series in \eqref{eq:behavPsi} as one wants, and this leads to arbitrarily large values for $c$. This means that the effects of the additional jumps are smaller than polynomial in $n$, which is basically a consequence of the fact that the jumps on the lens decrease exponentially as a function of $n$ away from the origin. 

For a $2\times 2$ RHP we have $D_\varphi = \operatorname{diag}[\varphi,-\varphi]$, where $\varphi$ is a function related to the normalization and opening of lenses. Usually, one can show that $\varphi^b$ defines a conformal map that maps $0$ to itself. Then we can construct a conformal map $f$, such that $n D_\varphi(z)+\theta(n^b f(z)) = 0$ exactly. Thus in the $2\times 2$ case \eqref{eq:defmathE} turns into $E(z) = N(z) D(z)^{-1} B(n^b f(z))^{-1}$. Under the assumption that $B$ has power law behavior, we infer that there exists a $d>0$ such that $E(z)=\mathcal O(n^\frac{d}{2})$ as $n\to\infty$, on any circle around the origin of fixed radius $r>0$. In such a case we have, taking $k=0$ in \eqref{eq:matchingAsympS}, that on $\partial D(0,r)$ as $n\to\infty$
\begin{align} \label{eq:concludingMatch}
E(z) \mathring P(z) N(z)^{-1} = E(z) \left(\mathbb I+\mathcal O(n^{-c})\right) E(z)^{-1} = \mathbb I + \mathcal O(n^{d-c}). 
\end{align}
It turns out that we generally have $c=b>d$, and we thus obtain the matching on $\partial D(0,r)$. In conclusion, for $2\times 2$ RHPs, $E$ really \textit{is} a sufficient prefactor to obtain the matching. In fact, this is why we added the inverse in \eqref{eq:matchingAsympS}. 

The case of an $m\times m$ RHP with $m>2$ is fundamentally different. It is generally not possible to arrange that $n D_\varphi(z)+\theta(n^b f(z)) = 0$ exactly. Frequently, the best we can get is that $n D_\varphi(z)+\theta(n^b f(z)) = \mathcal O(n z^\frac{1}{a})$ for some $a>0$, which originates from an expansion for the $\varphi$-functions in fractional powers of $z$. This means that, on a circle of fixed radius, $E$ varies wildly with $n$, exponentially in fact.  In larger size RHPs we therefore really have to use a shrinking circle to match $n D_\varphi(z)$ and $\theta(n^b f(z))$ appropriately, where $n^{-a}$ seems to be an optimal choice for the shrinking radius. Using a shrinking circle essentially has the effect of lowering $c$ in the right-hand side of \eqref{eq:concludingMatch}, and then the inequality $c>d$ might not hold anymore, and the matching is not achieved.


\subsection{Main result}

The main result of this paper is that, rather than obtaining the matching on a circle $\partial D(0,r)$, we can actually obtain a matching on two circles instead of one. This double matching is sufficient if one is interested in the large $n$ behavior of the solution of the RHP, although an additional, but natural, assumption is necessary for a scaling limit at the special point itself (see Lemma \ref{thm:scalingb}). 

\begin{theorem}[\textbf{main result}] \label{lem:matching}
Let $\mathring P$ and $N$ be defined in a neighborhood of $\overline{D(0,r)}$ for some $r>0$. These are matrix-valued functions of size $m\times m$ that may vary with $n$. 
Let $a, b, c, d, e \geq 0$ satisfy
\begin{align} \label{eq:assumpabcde}
a\leq e < b\quad\quad\text{ and }\quad\quad  d<\min(b,c).
\end{align}
Suppose that uniformly for $z\in\partial D(0,n^{-a})$ as $n\to\infty$
\begin{align}\label{eq:almostMatching}
\mathring P(z) N(z)^{-1} E(z) &= \mathbb I + \frac{C(z)}{n^b z} + \mathcal O\left(n^{-c}\right),
\end{align}
where $C$ and $E$ are $m\times m$ functions in a neighborhood of $\overline{D(0,r)}$ that may vary with $n$, and
\begin{itemize}
\item[(i)] $C$ is meromorphic with only a possible pole at $z=0$, whose order is bounded by some non-negative integer $p$ for all $n$, and $C$ is uniformly bounded for $z\in \partial D(0,n^{-a})$ as $n\to\infty$,
\item[(ii)] $E$ is non-singular, analytic, and uniformly for $z,w\in \partial D(0,n^{-a})$ we have as $n\to\infty$
\begin{align} \label{eq:assumpE}
E(z) = \mathcal O(n^\frac{d}{2}), \quad\quad E(z)^{-1} = \mathcal O(n^{\frac{d}{2}}),\quad\quad \text{ and }\quad E(z)^{-1} E(w) &= \mathbb I+\mathcal O(n^e (z-w)).
\end{align}
\end{itemize}
Then there are non-singular analytic functions  ${E_n^0:\overline{D(0,n^{-a})}\to\mathbb C^{m\times m}}$, ${E_n^\infty:\overline{A(0;n^{-a},\infty)}\to\mathbb C^{m\times m}}$ such that as $n\to\infty$ 
\begin{align} \label{eq:matchingonna}
E_n^0(z) \mathring P(z) &= \left(\mathbb I + \mathcal O(n^{d-c})\right) E_n^\infty(z) N(z), &\text{uniformly for }z\in \partial D(0,n^{-a}),\\ \label{eq:matchingonr}
E_n^\infty(z) &= \mathbb I + \mathcal O(n^{d-b}), &\text{uniformly for }z\in \partial D(0,r).
\end{align}
\end{theorem}

In Theorem \ref{lem:matching} and in the rest of the paper it is tacitly assumed that $r>1$ when $a=0$, to assure that $\partial D(0,n^{-a})$ is contained in $D(0,r)$. Notice that no (direct) smoothness condition is imposed on the functions $\mathring P$ and $N$. In particular, they are allowed to have jumps on certain curves and $N$ may blow up around possible special points of the RHP for $S$. It is also not part of the conditions that $\mathring P$ is of the form \eqref{eq:defmathringP}. The theorem allows for the possibility that the global parametrix $N$ depends on $n$. Notice also that the exponents $d-c$ and $b-c$ in the statement are negative due to the assumptions. There is a workaround when the inequality $b>d$ is violated, see Remark \ref{remark:globalMatching}. 

We mention that the prefactors $E_n^0$ and $E_n^\infty$ in Theorem \ref{lem:matching} are not unique. Indeed, any multiplication of $E_n^0$ and $E_n^\infty$ by analytic functions with behaviors $\mathbb I + \mathcal O(n^{d-c})$ and $\mathbb I + \mathcal O(n^{d-b})$ respectively on the corresponding circles will again yield appropriate prefactors to obtain the matchings \eqref{eq:matchingonna} and \eqref{eq:matchingonr}. 

In the situation of \eqref{eq:matchingAsympS} we may put
\begin{align} \label{eq:naturalC}
C(z) = \frac{z}{f(z)}\left(C_1  +\frac{C_2}{(n^b z)} \frac{z}{f(z)} + \ldots \frac{C_k}{(n^b z)^{k-1}} \frac{z^{k-1}}{f(z)^{k-1}} \right),
\end{align}
and this function is indeed uniformly bounded for $z\in\partial D(0,n^{-a})$ when $a <b$. As mentioned before, we can take as many terms of \eqref{eq:behavPsi} as we want, and in doing so make $c$, which in general equals $(b-a) (k+1)$, as big as we want. Essentially, the assumption $a\leq e<b$ is really the only important one in \eqref{eq:assumpabcde}. 

When the assumptions of Theorem \ref{lem:matching} are met, one can define the final transformation $R$ as follows.
\begin{align} \label{finalT}
R(z) = \left\{\begin{array}{ll}
S(z) N(z)^{-1}, & z\in A(0;r,\infty),\\
S(z) N(z)^{-1} E_n^\infty(z)^{-1}, & z\in A(0;n^{-a},r),\\
S(z) \mathring P(z)^{-1} E_n^0(z)^{-1}, & z\in D(0,n^{-a}).
\end{array}\right.
\end{align}
Here again, for convenience we absorbed any local parametrices around non-zero special points into the definition of $N$. As one can check, $R$ will have jump matrices on the inner and outer circle of the form $\mathbb I+\mathcal O(n^{d-c})$ and $\mathbb I+\mathcal O(n^{d-b})$ respectively. In particular, we have a matching condition on two circles.

In a way, we are in the situation of \eqref{eq:finalOriginal} if we identify the local parametrix $P$ with
\begin{align} \nonumber
P(z) = \left\{\begin{array}{ll}
E_n^\infty(z) N(z), & z\in A(0;n^{-a},r),\\
E_n^0(z) \mathring P(z), & z\in D(0,n^{-a}).
\end{array}\right.
\end{align}
In this approach, I suspect that it will generally not be hard to find out what one should pick for $E$ in Theorem \ref{lem:matching}, for example by explicitly writing the asymptotic behavior of $\Psi$ in the form \eqref{eq:behavPsi} and then imposing \eqref{eq:defmathE}. The main technical part to apply the theorem will be to prove that the estimates in \eqref{eq:assumpE} are satisfied. In Section \ref{sec:examples} we will show for several examples in the literature that one obtains the situation \eqref{eq:matchingAsympS}, and we will show how to identify $E$ and prove that the estimates in \eqref{eq:assumpE} hold.

\section{Proof of the main result} \label{sec:proofMainResult}

In this section we always assume that the conditions of Theorem \ref{lem:matching} are met. In particular, we will use $\mathring P, N, E, C$ and $r,a,b,c,d,e, p$ without reference, and they will satisfy the conditions of Theorem \ref{lem:matching}. To avoid confusion, we mention that it is not assumed that $C$ is necessarily of the form \eqref{eq:naturalC}. 

It is not hard to see that, when $c\leq b-a$, prefactors that satisfy the properties of Theorem \ref{lem:matching} are given by $E_n^0=E$ and $E_n^\infty = \mathbb I$. In other words, a double matching is unnecessary in such a case. Henceforth, we will exclude this trivial case. 

\begin{assumption} \label{cbiggera-b}
We have $c>b-a$. 
\end{assumption}

\subsection{Definition of the analytic prefactors} \label{sec:defPrefactors}

In this subsection we give a definition of the analytic prefactors $E_n^0$ and $E_n^\infty$. In the remaining subsections we will prove that they indeed satisfy the desired properties from Theorem \ref{lem:matching}. 

Suppose that $F$ is a (matrix valued) function that has a Laurent series expansion around $z=0$, i.e.,
\begin{align*}
F(z) = \sum_{k=-\infty}^\infty F_k z^k
\end{align*}
for certain coefficients $F_k$ in a punctured disc around $z=0$. We denote by $F^-$ the principal part of $F$, 
\begin{align*}
F^-(z) = \sum_{k=-\infty}^{-1} F_k z^k.
\end{align*}
By $F^+$ we denote the regular part of $F$, that is, $F^+ = F - F^-$. The regular and principal part can neatly be expressed using a Cauchy-operator. Namely, if $F$ is analytic on some open neighborhood of $D(0,\rho)\setminus\{0\}$ that does not contain $0$, for some $\rho>0$, then we have
\begin{align} \nonumber
F^+(z) &= \frac{1}{2\pi i} \oint_{\partial D(0,\rho)} \frac{F(s)}{s-z} ds, & |z|<\rho,\\ \label{eq:principalPartInt}
F^-(z) &= -\frac{1}{2\pi i} \oint_{\partial D(0,\rho)} \frac{F(s)}{s-z} ds, & |z|>\rho,
\end{align}
on this neighborhood. In fact, if $F$ has a pole of order at most $q$ in $z=0$, then we deduce that
\begin{align} \label{eq:sumCauchyOp}
F^-(z) &= \frac{1}{2\pi i} \sum_{j=1}^q \left(\oint_{\partial D(0,\rho)} F(s) s^j \frac{ds}{s}\right) z^{-j}
\end{align}
for all $z\neq 0$. We will limit ourselves to the case of meromorphic $m\times m$ functions $F$ on $\overline{D(0,r)}$ that have at most a pole in $z=0$ and no other singularities. For such functions we define an operator $\pi$ through
\begin{align} \label{def:operator}
\pi F = - F^+ F - F F^- + F^+ F^- + F^+ F F^-.
\end{align}
Notice that $\pi F$ is again a meromorphic function on $\overline{D(0,r)}$ that has at most a pole in $z=0$ and no other singularities. We mention that $\pi F$ has a pole of order at most $2 q$ if $F$ has a pole of order $q$. 

The meromorphic function that will be of interest to us is defined as follows. 

\begin{definition} \label{def:F}
On $\overline{D(0,r)}$ we define the $n$-dependent function
\begin{align} \label{eq:defF}
F(z) = \displaystyle \frac{E(z) C(z) E(z)^{-1}}{n^b z}.
\end{align} 
\end{definition}

We emphasize that $F$ is a meromorphic function with only a possible pole at $z=0$, whose order is bounded by $p+1$. Notice that the $n$-dependence of $F$ is not only due to the $n^b$ factor, but also due to $E$ and $C$. 
This function $F$ originates from the following. 

\begin{proposition} \label{prop:0step} 
We have uniformly for $z\in\partial D(0,n^{-a})$ that as $n\to\infty$
\begin{align} \label{eq:0stepF}
E(z) \mathring P(z) N(z)^{-1} = \mathbb I + F(z) + E(z) \mathcal O(n^{-c}) E(z)^{-1},
\end{align}
where $F$ is as in Definition \ref{def:F}.
\end{proposition}

\begin{proof}
We can left-multiply the entire equation \eqref{eq:almostMatching} with $E(z)$ and right-multiply it with $E(z)^{-1}$. Comparing with Definition \ref{def:F} then immediately yields the result. 
\end{proof}

Notice that the first two  estimates in \eqref{eq:assumpE} of Theorem \ref{lem:matching} imply that $E(z) \mathcal O(n^{-c}) E(z)^{-1}=\mathcal O(n^{d-c})$ uniformly for $z\in\partial D(0,n^{-a})$ as $n\to\infty$. Since $d<c$ by the assumptions of Theorem \ref{lem:matching}, such an expression is thus small. 

The idea now is as follows. Starting with \eqref{eq:0stepF}, we left-multiply by $\mathbb I - F^+(z)$ and right-multiply by $\mathbb I - F^-(z)$, to give after some bookkeeping that
\begin{multline}
\left(\mathbb I - F^+(z)\right) E(z) \mathring P(z) N(z)^{-1} \left(\mathbb I - F^-(z)\right) = \mathbb I + \pi F(z)\\
+ \left(\mathbb I - F^+(z)\right) E(z) \mathcal O(n^{-c}) E(z)^{-1} \left(\mathbb I - F^-(z)\right),
\end{multline}
where the term in the last line turns out to be $\mathcal O(n^{d-c})$ on $\partial D(0,n^{-a})$ as $n\to\infty$. Repeating the iteration, we get
\begin{multline}
\left(\mathbb I - (\pi F)^+(z)\right) \left(\mathbb I - F^+(z)\right) E(z) \mathring P(z) N(z)^{-1} \left(\mathbb I - F^-(z)\right) \left(\mathbb I - (\pi F)^-(z)\right) = \mathbb I + \pi^2 F(z)\\
+ \left(\mathbb I - (\pi F)^+(z)\right) \left(\mathbb I - F^+(z)\right) E(z) \mathcal O(n^{-c}) E(z)^{-1} \left(\mathbb I - F^-(z)\right) \left(\mathbb I - (\pi F)^-(z)\right),
\end{multline}
where again the term in the last line turns out to be $\mathcal O(n^{d-c})$. After $k$ steps we will obtain $\mathbb I + \pi^k F$ and a remainder. Our goal will be to prove that these remainders are indeed $\mathcal O(n^{d-c})$ on $\partial D(0,n^{-a})$, and that $\pi^k F$ will be small enough on $\partial D(0,n^{-a})$ for sufficiently big (but fixed) $k$ as $n\to\infty$.

Asuming that this reasoning is correct, it justifies the following explicit definition for prefactors $E_n^0$ and $E_n^\infty$ that will satisfy the required properties for Theorem \ref{lem:matching}, although an argument for the matching on the outer circle is still missing. As it turns out, the matching on the outer circle will follow without much trouble. 

\begin{definition} \label{def:prefactors}
Let $\pi$ be as in \eqref{def:operator} and let $F$ be as in \eqref{eq:defF}. Let $K$ be the biggest integer such that 
\begin{align*} 
2^{K} < \frac{a+c-e}{b-e}.
\end{align*}
We define 
\begin{align} \label{eq:defEn+}
E_n^0(z) &= \prod_{j=0}^K \left(\mathbb I - (\pi^{K-j} F)^+(z)\right) E(z), & z\in \overline{D(0,n^{-a})},\\ \label{eq:defEn-}
E_n^\infty(z) &= \left(\prod_{j=0}^K \left(\mathbb I - (\pi^j F)^-(z)\right)\right)^{-1}, & z\in \overline{A(0;n^{-a},\infty)}.
\end{align}
\end{definition}

For clarity, the notation of the products in \eqref{eq:defEn+} and \eqref{eq:defEn-} means that a factor with $j=0$ is at the left and a factor with $j=K$ is at the right. We remark that $K$ is well-defined and non-negative due to Assumption \ref{cbiggera-b}. Notice that $E_n^\infty(z)^{-1}$ is in fact a polynomial evaluated in $1/z$ with constant term $\mathbb I$. This polynomial has degree at most $2^{\frac{K(K+1)}{2}} (p+1)$. We shall prove later that $E_n^0$ and $E_n^\infty$ are indeed non-singular analytic functions, provided that $n$ is big enough (see Proposition \ref{prop:En0inftyWelld}). In particular, we show that $E_n^\infty$ is well-defined, i.e., that we may take the inverse of the product in \eqref{eq:defEn-} when $n$ is big enough.

\subsection{The structure of iterations of $\pi$}

It is instructive to write down explicitly what the effect of $\pi$ on $F$ is. We write $q=p+1$. First, we notice that we can eliminate $F^+$ in \eqref{def:operator}. Namely, using $F^+ = F - F^-$, we can rewrite
\begin{align*}
\pi F = -F^2 + F^- F - (F^-)^2 + F^2 F^- - F^- F F^-.
\end{align*}
We may use \eqref{eq:sumCauchyOp} with any $0<\rho\leq r$ to write this as
\begin{align*} 
(\pi F)(z) =& - F(z)^2\\ 
&+ \sum_{j_1=1}^q \frac{1}{2\pi i} \oint_{\partial D(0,\rho)}  F(s_1) F(z) \left(\frac{s_1}{z}\right)^{j_1} \frac{ds_1}{s_1} \\ 
&- \sum_{j_1=1}^q \sum_{j_2=2}^q \frac{1}{(2\pi i)^2} \oint_{\partial D(0,\rho)^2} F(s_1) F(s_2) \left(\frac{s_1}{z}\right)^{j_1} \left(\frac{s_2}{z}\right)^{j_2} \frac{ds_1}{s_1} \frac{ds_2}{s_2}\\ 
&+ \sum_{j_1=1}^q \frac{1}{2\pi i} \oint_{\partial D(0,\rho)}  F(z)^2 F(s_1) \left(\frac{s_1}{z}\right)^{j_1}\frac{ds_1}{s_1}\\ 
&+ \sum_{j_1=1}^q \sum_{j_2=2}^q\frac{1}{(2\pi i)^2} \oint_{\partial D(0,\rho)^2} F(s_1) F(z) F(s_2) \left(\frac{s_1}{z}\right)^{j_1} \left(\frac{s_2}{z}\right)^{j_2} \frac{ds_1}{s_1} \frac{ds_2}{s_2}.
\end{align*}
So we have five different expressions, four of which are accompanied by integrals and sums. If we would write down the next iteration (which is a cumbersome task) the number of integrals and summations accompanying such expressions may be bigger, and the summations may also go from $1$ to $2q$ or even $3q$, rather than $q$, since $\pi$ may increase the order of the poles. 
A general structure for $\pi^k F$ (that follows from repeated application of Proposition \ref{prop:OFkGH}(c) in the next subsection) emerges. Namely, using the notation $s_0=z$, any $\pi^k F$ is a linear combination of expressions of the form
\begin{align} \label{eq:piFsumint4}
\frac{1}{(2\pi i)^t} \int_{\partial D(0,\rho)^t} F\left(s_{i_1}\right) \cdots F\left(s_{i_l}\right) \left(\frac{s_1}{z}\right)^{j_1} \cdots \left(\frac{s_t}{z}\right)^{j_t} \frac{ds_1}{s_1} \cdots \frac{ds_t}{s_t},
\end{align}
for some positive integer $2^k\leq l\leq 3^k$, some integer $0\leq t\leq l$, integers $i_1,\ldots,i_l\in\{0,1,\ldots,t\}$, and integers $j_1,\ldots,j_t$ that are not necessarily positive for further iterations.  The number of expressions of the form \eqref{eq:piFsumint4} in $\pi^k F$ grows double exponentially as a function of $k$, and is bounded by an $n$-independent constant for fixed $k$ in particular. We may in fact take $t = l$ by artificially adding integrations. It is thus equivalent to replace \eqref{eq:piFsumint4} by the somewhat more appealing expression
\begin{align} \label{eq:piFsumint2}
\frac{1}{(2\pi i)^l} \oint\cdots\oint F\left(s_{i_1}\right) \cdots F\left(s_{i_l}\right) \left(\frac{s_1}{z}\right)^{j_1} \cdots \left(\frac{s_l}{z}\right)^{j_l}  \frac{ds_1}{s_1} \cdots \frac{ds_l}{s_l},
\end{align}
where the integrations are over any circle in $\overline{D(0,r)}$ around $z=0$, and now $i_1,\ldots,i_l\in\{0,1,\ldots,l\}$. We choose to leave out the factors $(2\pi i)^{-l}$ from now on, as they will only be extra baggage in what follows. Since these expressions will play a key role, we introduce a notation for them.

\begin{definition} \label{def:Ln}
Let $F$ be as in Definition \ref{def:F} and let $l$ be a positive integer. Write $s_0=z$.\\
For $i=(i_1,\ldots,i_l)\in \{0,1,\ldots,l\}^l$ and $j=(j_1,\ldots,j_l)\in \mathbb Z^l$, we define the function
\begin{align} \label{eq:piFsumint3}
I^{[l]}_{i,j}(n, z) =  \oint\cdots\oint F\left(s_{i_1}\right) \cdots F\left(s_{i_l}\right) \left(\frac{s_1}{z}\right)^{j_1} \cdots \left(\frac{s_l}{z}\right)^{j_l}  \frac{ds_1}{s_1} \cdots \frac{ds_l}{s_l},
\end{align} 
where $n=1,2,\ldots$ and $z\in\overline{D(0,r)}$. The integrations are over any circle in $\overline{D(0,r)}$ around the origin, with positive orientation.
\end{definition}

We emphasize that the $n$-dependence in \eqref{eq:piFsumint3} comes from $F$, as is clear from Definition \ref{def:F}. So far, we have not used any information about $F$, except for the order of its pole. To understand the behavior of the $I^{[l]}_{i,j}(n,z)$ on $\partial D(0,n^{-a})$ as $n\to\infty$, we need to understand the behavior of products of the form $F\left(s_{i_1}\right) \cdots F\left(s_{i_l}\right)$ on $\partial D(0,n^{-a})^l$. It will turn out that in any linear combination of expressions of the form \eqref{eq:piFsumint2}, the terms with the smallest number of $F$ factors are dominant.

\begin{lemma} \label{prop:behavA}
Let $l=1,2,\ldots$ and let $F$ be as in Definition \ref{def:F}. We have
\begin{align}
F(z_1) \cdots F(z_l) = \mathcal O(n^{a+d-e-(b-e)l}),
\end{align}
uniformly for $z_1,\ldots,z_l \in \partial D(0,n^{-a})$ as $n\to\infty$.
\end{lemma}

\begin{proof}
(a) By the assumptions of Theorem \ref{lem:matching} we have uniformly for $z\in \partial D(0,n^{-a})$ that as $n\to\infty$
\begin{align} \label{eq:EE-C}
E(z) = \mathcal O(n^\frac{d}{2}), \quad\quad E(z)^{-1} = \mathcal O(n^\frac{d}{2}), \quad\quad C(z) = \mathcal O(1). 
\end{align}
Now using \eqref{eq:defF}, we have uniformly for $z\in \partial D(0,n^{-a})$ as $n\to\infty$
\begin{align} \nonumber
F(z) &= \mathcal O(n^{a-b+\frac{d}{2}+\frac{d}{2}}) = \mathcal O(n^{a+d-e-(b-e)}).
\end{align}
This proves the statement for $l=1$. By the assumptions of Theorem \ref{lem:matching} we also have uniformly for $z,w\in \partial D(0,n^{-a})$ as $n\to\infty$
\begin{align} \label{eq:E-Ee-a}
E(z)^{-1} E(w) = \mathcal O(n^{e-a}).
\end{align}
Combining \eqref{eq:E-Ee-a} with the boundedness of $C$ we have uniformly for $z_1,\ldots,z_l\in \partial D(0,n^{-a})$ as $n\to\infty$
\begin{multline} \label{eq:proofAA}
F(z_1) \cdots F(z_l) = (n^b z)^{-l} E(z_1) C(z_1)\\ \left(E(z_1)^{-1} E(z_2)\right) C(z_2)  \cdots 
C(z_{l-1}) \left(E(z_{l-1})^{-1} E(z_l)\right) C(z_l) E(z_l)^{-1}
\end{multline}
\vspace{-0.5cm}
\begin{align} \nonumber
\hspace{2.3cm}= \mathcal O(n^{(a-b)l+\frac{d}{2}+(e-a)(l-1)+\frac{d}{2}}) 
= \mathcal O(n^{a+d-e-(b-e) l}).
\end{align}
\end{proof}

\begin{corollary} \label{cor:growthInt}
Let $l$ be a positive integer, let $i\in \{0,1,\ldots,l\}^l$ and let $j\in \mathbb Z^l$. Then we have
\begin{align}
I^{[l]}_{i,j}(n,z) = \mathcal O\left(n^{a+d-e-(b-e)l}\right),
\end{align}
uniformly for $z\in \partial D(0,n^{-a})$ as $n\to\infty$, where $I^{[l]}_{i,j}$ is as in Definition \ref{def:Ln}. 
\end{corollary}

\begin{proof}
We may take the integrations in \eqref{eq:piFsumint3} over $\partial D(0,n^{-a})$. Indeed, for $s_1,\ldots,s_l\in\partial D(0,n^{-a})$,
\begin{align}
\left(\frac{s_1}{z}\right)^{j_1}, \ldots, \left(\frac{s_l}{z}\right)^{j_l} = \mathcal O(1)
\quad\quad\text{ and }\quad\quad \frac{ds_1}{s_1}, \ldots, \frac{ds_l}{s_l} = \mathcal O(1). 
\end{align}
Combining these with Lemma \ref{prop:behavA} yields the result.
\end{proof}

Since $b>e$ it follows from Corollary \ref{cor:growthInt} that the dominant behavior of any iteration $\pi^k F$ is given by the expressions of the form \eqref{eq:piFsumint3} that have the smallest number of $F$ factors. 

\subsection{Definition of auxiliary spaces}

As stated in Subsection \ref{sec:defPrefactors}, there are two things that we should prove in order to obtain the matching on the inner circle. Namely, we want to prove that $\pi^K F$ (with $K$ as in Definition \ref{def:prefactors}) is small on $\partial D(0,n^{-a})$, and we want to prove that the remainder is small on $\partial D(0,n^{-a})$. To be more accurate, we want the following. 
\begin{itemize}
\item[(i)] Uniformly on $\partial D(0,n^{-a})$ as $n\to\infty$, we have
\begin{align} \label{eq:piKFTo0}
\pi^K F(z) = \mathcal O(n^{d-c}).
\end{align}
\item[(ii)] Uniformly on $\partial D(0,n^{-a})$ as $n\to\infty$, we have
\begin{align} \nonumber
\left(\mathbb I - (\pi^K F)^+(z)\right) &\cdots \left(\mathbb I - F^+(z)\right) E(z) \mathcal O(n^{-c}) E(z)^{-1} \left(\mathbb I - F^-(z)\right) \cdots \left(\mathbb I - (\pi^K F)^-(z)\right)\\ \label{eq:remainderTo0}
&= \mathcal O(n^{d-c}).
\end{align}
Here we understand that the function implied by the $\mathcal O(n^{-c})$ term is the same as the function implied by the $\mathcal O(n^{-c})$ term in \eqref{eq:almostMatching}, although its specifics turn out to be irrelevant. 
\end{itemize}

For any positive integer $k$, $\pi^k F$ is a linear combination of the $I^{[l]}_{i,j}$ as in Definition \ref{def:Ln} (again, this will follow from repeated application of Proposition \ref{prop:OFkGH}(c), that we prove in a moment), where $l$ ranges between $2^k$ and $3^k$. The coefficients in the linear combinations are $\pm (2\pi i)^{-l}$. In principle, we can try to be as precise as possible about the structure of $\pi^k F$ when proving \eqref{eq:piKFTo0}, but there is a good reason not to. Namely, to prove \eqref{eq:remainderTo0}, we need to understand expressions such as
\begin{align*}
\left(\mathbb I - (\pi^K F)^+(z)\right) &\cdots \left(\mathbb I - F^+(z)\right) \quad\text{ and }\quad
\left(\mathbb I - F^-(z)\right) \cdots \left(\mathbb I - (\pi^K F)^-(z)\right).
\end{align*}
Hence, besides $\pi^k F, (\pi^k F)^+$ and $(\pi^k F)^-$ with $k=0,1,\ldots,K$, we should also understand sums and products of these. There is a natural framework to accomplish just that, which is to imbed the functions of interest in a conveniently chosen ring. This ring will be $\mathcal O_F^0$, as defined below.

\begin{definition} \label{def:OFk}
We define $O_F^0 = \{\lambda \mathbb I \hspace{0.2cm}|\hspace{0.2cm} \lambda\in\mathbb C\}$. For any positive integer $l$, we define $O_F^l$ as the linear span of $\{I^{[l]}_{i,j} \hspace{0.2cm}|\hspace{0.2cm} i\in \{0,1,\ldots,l\}^l, j\in\mathbb Z^l\}$, where $I^{[l]}_{i,j}$ is as in Definition \ref{def:Ln}.\\
For any non-negative integer $l$ we define $\mathcal O_F^l$ as the linear span of $O_F^l \cup O_F^{l+1} \cup O_F^{l+2} \cup \ldots$
\end{definition}


We remind the reader that the linear span consists of all finite linear combinations. It is clear from the definition that $\ldots \subset \mathcal O_F^2 \subset \mathcal O_F^1 \subset \mathcal O_F^0$. Notice that $\mathcal O_F^k$ is a vector space by construction. We deduce from the next proposition that $\mathcal O_F^k$ is in fact a graded ring, which is unital when $k=0$. 

\begin{proposition} \label{prop:OFkGH}
Let $k,l=0,1,\ldots$ and let $G\in\mathcal O_F^k$ and $H\in\mathcal O_F^l$. We have
\begin{itemize}
\item[(a)] $G H \in \mathcal O_F^{k+l}$.
\item[(b)] $G^+, G^-\in \mathcal O_F^k$.
\item[(c)] $\pi G\in\mathcal O_F^{2k}$.
\item[(d)] $G(n,z) = \mathcal O\left(n^{a+d-e-(b-e)k}\right)$ when $k\geq 1$, uniformly for $z\in\partial D(0,n^{-a})$ as $n\to\infty$.
\item[(e)] $G(n,z) = \mathcal O(n^{d-c})$ when $k\geq \frac{a+c-e}{b-e}$, uniformly for $z\in\partial D(0,n^{-a})$ as $n\to\infty$.
\end{itemize}
\end{proposition}

\begin{proof} We omit the cases $k=0$ or $l=0$, as these are similar to the cases $k,l\geq 1$, but easier. \\
(a) A product of two expressions of the form \eqref{eq:piFsumint3} is again such an expression (with the indices relabeled), but the number of $F$ factors is now the sum of the number of $F$ factors of each of the two factors.\\
(b) It suffices to prove that $G^-\in\mathcal O_F^k$. Without loss of generality we may assume that $G = I^{[l]}_{i,j}$ for some choice of $i$ and $j$, that is,
\begin{align*}
G(n,z) = \oint\cdots\oint  F\left(s_{i_1}\right) \cdots F\left(s_{i_l}\right) \left(\frac{s_1}{z}\right)^{j_1} \cdots \left(\frac{s_l}{z}\right)^{j_l}  \frac{ds_1}{s_1} \cdots \frac{ds_l}{s_l},
\end{align*}
with all expressions as in Definition \ref{def:Ln}, and with $l\geq k$. For each $n$, $G$ has a pole of order at most $q = j_1+\ldots+j_l+(p+1) n_i$, where $n_i$ is the number of $0$ components of $i$. If $q\leq 0$, then $G$ is analytic and we have $G^- = 0$. In that case we are done. So let us assume that $q\geq 1$ henceforth. 

If $(i_1,\ldots,i_l)$ is a permutation of $(1,2,\ldots,l)$, then $G$ is simply a multiple of $z^{-q}$. This would imply that $G^-=G$, and we are done in that case. Let us thus assume that at least one of the integrations is artificial, i.e., let us assume that we can eliminate the integration with respect to, say, $ds_{l}$. Since we have assumed that $G$ has a pole, $G$ is in particular not identically $0$ and we then necessarily have $j_l = 0$. The corresponding integral thus yields a factor $2\pi i$. Now using \eqref{eq:sumCauchyOp}, where we redefine $s_l=s$, we have
\begin{align*}
G^-(n,z) = \sum_{j=1}^q \oint\cdots\oint  F\left(s_{i_1'}\right) \cdots F\left(s_{i_l'}\right) \left(\frac{s_1}{s_l}\right)^{j_1} \cdots \left(\frac{s_{l-1}}{s_l}\right)^{j_{l-1}} \left(\frac{s_l}{z}\right)^{j} \frac{ds_1}{s_1} \cdots \frac{ds_{l-1}}{s_{l-1}} \frac{ds_l}{s_l},
\end{align*}
where $i_t' = l$ when $i_t= 0$ and $i_t'=i_t$ otherwise. Trivially, we have
\begin{align*}
\left(\frac{s_1}{s_l}\right)^{j_1} \cdots \left(\frac{s_{l-1}}{s_l}\right)^{j_{l-1}} \left(\frac{s_l}{z}\right)^{j} = \left(\frac{s_1}{z}\right)^{j_1} \cdots \left(\frac{s_l}{z}\right)^{j_l} 
\left(\frac{s_l}{z}\right)^{j-(j_1+\ldots+j_{l-1})}.
\end{align*}
We conclude that $$G^- = \sum_{j=1}^q I^{[l]}_{i',(j_1,\ldots,j_{l-1}, j-(j_1+\ldots+j_{l-1}))}\in\mathcal O_F^k.$$\\
(c) Starting from the definition \eqref{def:operator}, this is a direct consequence of (a) and (b).\\
(d) This is a direct consequence of Corollary \ref{cor:growthInt}.\\
(e) This follows from (d) and the assumption $b>e$. 
\end{proof}

To also treat the remainder terms that we discussed in Subsection \ref{sec:defPrefactors}, we have to define yet another space. 
The intuition behind this space, is that we want it to be invariant under multiplication with $\mathbb I - (\pi^k F)^\pm$ for any $k$. This property should be satisfied by the remainders (when multiplied from the correct side). In the definition that follows, we will denote by $O$ the function on $\overline{D(0,r)}$ implied by the $\mathcal O(n^{-c})$ term in \eqref{eq:almostMatching}, i.e.,
\begin{align} \nonumber
O(z) = \mathring P(z) N(z)^{-1} E(z) - \mathbb I - \frac{C(z)}{n^b z}, \hspace{3cm} z\in\overline{D(0,r)}.
\end{align}
We will not need the explicit description for $O(z)$ though. 

\begin{definition}
Let $O$ be as above. We define $\mathcal O_E$ as the linear span of $\{G E O E^{-1} H \hspace{0.2cm}|\hspace{0.2cm} G, H \in \mathcal O_F^0 \}$, with $\mathcal O_F^0$ as in Definition \ref{def:OFk}.
\end{definition}

Before proving some properties of $\mathcal O_E$, we will need the following intermediate result. 

\begin{proposition} \label{prop:behavA2}
Let $F$ be as in Definition \ref{def:F}. We have 
\begin{align} \nonumber
F(z_1) E(z_2) = E(z_1)\mathcal O(n^{e-b})\quad\quad\text{and}\quad\quad
E(z_1)^{-1} F(z_2) = \mathcal O(n^{e-b}) E(z_2)^{-1},
\end{align}
uniformly for $z_1,z_2\in\partial D(0,n^{-a})$ as $n\to\infty$.
\end{proposition}

\begin{proof}
This follows essentially from the same arguments that we used to prove Lemma \ref{prop:behavA}. 
\end{proof}

\begin{proposition} \label{prop:prodOAE}
Let $G \in \mathcal O_E$ and let $H\in\mathcal O_F^0$. We have
\begin{itemize}
\item[(a)] $G(n,z) = \mathcal O(n^{d-c})$ uniformly for $z\in\partial D(0,n^{-a})$ as $n\to\infty$.
\item[(b)] $H G \in \mathcal O_E$ and $G H \in \mathcal O_E$. 
\end{itemize}
\end{proposition}

\begin{proof}
(a) We may assume without loss of generality that
\begin{multline*}
G(n,z) = \oint\cdots\oint F\left(s_{i_1}\right) \cdots F\left(s_{i_k)}\right)  \left(\frac{s_1}{z}\right)^{j_1} \cdots \left(\frac{s_k}{z}\right)^{j_k}  \frac{ds_1}{s_1} \cdots \frac{ds_k}{s_k}
E(z) O(z) E(z)^{-1} \\
\oint\cdots\oint F\left(s_{i_{k+1}}\right) \cdots F\left(s_{i_{k+l}}\right)  \left(\frac{s_{k+1}}{z}\right)^{j_{k+1}} \cdots \left(\frac{s_{k+l}}{z}\right)^{j_{k+l}}  \frac{ds_{k+1}}{s_{k+1}} \cdots \frac{ds_{k+l}}{s_{k+l}}
\end{multline*}
for some positive integers $k,l\geq 1$, where $F$ is as in Definition \ref{def:F}, $i_1,\ldots,i_{k}\in\{0,1,\ldots,k\}$ and $i_{k+1},\ldots,i_{k+l}\in \{0,k+1,k+2,\ldots,k+l\}$. As usual, we put $s_0=z$ and the integrations may be over any circle in $\overline{D(0,r)}$ around $0$. Repeated application of Proposition \ref{prop:behavA2} yields that uniformly for $z,z_1,\ldots,z_l\in \partial D(0,n^{-a})$ as $n\to\infty$
\begin{align*}
F(z_1) \cdots F(z_k) E(z) &= F(z_1) \cdots F(z_{k-1}) E(z_k) \mathcal O\left(n^{e-b}\right)\\
&= F(z_1) \cdots F(z_{k-2}) E(z_{k-1}) \mathcal O\left(n^{2(e-b)}\right)
= \ldots = E(z_1) \mathcal O\left(n^{k(e-b)}\right).
\end{align*}
Analogously, we have uniformly for $z,z_{k+1},\ldots,z_{k+1+l}\in \partial D(0,n^{-a})$ as $n\to\infty$
\begin{align*}
E(z)^{-1} F(z_{k+1}) \cdots F(z_{k+l}) = \mathcal O\left(n^{l(e-b)}\right) E(z_{k+l})^{-1}.
\end{align*}
Hence, taking all integrations over $\partial D(0,n^{-a})$, we have uniformly for $z\in\partial D(0,n^{-a})$ as $n\to\infty$
\begin{align} \label{eq:ointO}
G(n,z) = \oint \cdots \oint E\left(s_{i_1}\right) \mathcal O\left(n^{(k+l)(e-b)-c}\right) E\left(s_{i_{k+l}}\right) \left(\frac{s_1}{z}\right)^{j_1} \cdots \left(\frac{s_{k+l}}{z}\right)^{j_{k+1+l}}
\frac{ds_1}{s_1} \cdots \frac{ds_{k+l}}{s_{k+l}}.
\end{align}
Here we have used that $O(z)=\mathcal O(n^{-c})$ uniformly on $\partial D(0,n^{-a})$. From the assumptions of Theorem \ref{lem:matching} we have that
\begin{align} \label{eq:Ez1Ezk1l}
E\left(s_{i_1}\right) = \mathcal O(n^\frac{d}{2})\quad\quad\text{and}\quad\quad E\left(s_{i_{k+l}}\right)^{-1} = \mathcal O(n^\frac{d}{2})
\end{align}
as $n\to\infty$, uniformly for $s_{i_{1}},s_{i_{k+l}}\in\partial D(0,n^{-a})$. We also make the observation that
\begin{align*}
\left(\frac{s_1}{z}\right)^{j_1}, \ldots, \left(\frac{s_{k+l}}{z}\right)^{j_{k+l}} = \mathcal O(1)
\quad\quad\text{and}\quad\quad
\frac{ds_1}{s_1}, \ldots, \frac{ds_{k+l}}{s_{k+l}} = \mathcal O(1)
\end{align*}
uniformly on $\partial D(0,n^{-a})$. Combining this with \eqref{eq:Ez1Ezk1l} and \eqref{eq:ointO} yields $G(n,z) = \mathcal O\left(n^{d-c+(k+l)(e-b)}\right)$ uniformly for $z\in\partial D(0,n^{-a})$ as $n\to\infty$. Since $e<b$ by the assumptions of Theorem \ref{lem:matching}, it follows that $G(n,z) = \mathcal O(n^{d-c})$ uniformly for $z\in\partial D(0,n^{-a})$ as $n\to\infty$. The cases where $E(z) O(z) E(z)^{-1}$ is multiplied on the left or the right with $\mathbb I$ rather than a term of the form \eqref{eq:piFsumint3} are analogous, but easier.\\
(b) This is a direct consequence of Proposition \ref{prop:OFkGH}(a) and the fact that $\mathcal O_F^l\subset \mathcal O_F^0$ for all $l\geq 0$. 
\end{proof}

We note that Proposition \ref{prop:prodOAE}(b) implies that $\mathcal O_E$ is an $\mathcal O_F^0$-bimodule. Therefore, $\mathcal O_E$ will have the interpretation of a ``garbage can'', if you will, in what remains of the proof of Theorem \ref{lem:matching}. 

\begin{proposition}
Let $k$ be a positive integer and let $G\in\mathcal O_F^k$. Then $\pi G\in \mathcal O_F^{2 k}$ and 
\begin{align*}
(\mathbb I - G^+) \left(\mathbb I + G + \mathcal O_E\right) (\mathbb I - G^-)
= \mathbb I + \pi G + \mathcal O_E.
\end{align*}
\end{proposition}

\begin{proof}
Some straightforward bookkeeping yields
\begin{align*}
(\mathbb I - G^+) \left(\mathbb I + G + \mathcal O_E\right) (\mathbb I - G^-)
=& \mathbb I - G^+ - G^- + G^+ G^- + G - G^+ G - G G^- + G^+ G G^-\\
&+ (\mathbb I - G^+) \mathcal O_E (\mathbb I - G^-).
\end{align*}
By Proposition \ref{prop:prodOAE}(b) the term in the last line is $\mathcal O_E$. Using $G = G^-+G^+$ and the definition \eqref{def:operator} of $\pi$,  the proposition follows. 
\end{proof}

\begin{corollary} \label{cor:induction}
Let $k$ be a positive integer and let $G\in\mathcal O_F^1$. Then $\pi^{k} G\in \mathcal O_F^{2^k}$ and
\begin{align*}
\left(\prod_{j=0}^{k-1} \left(\mathbb I - (\pi^{k-j} G)^+\right)\right) \left(\mathbb I + G + \mathcal O_E\right) \left(\prod_{j=0}^{k-1} \left(\mathbb I - (\pi^{j} G)^-\right)\right)
= \mathbb I + \pi^{k} G + \mathcal O_E.
\end{align*}
\end{corollary}

\subsection{Proof of Theorem \ref{lem:matching}}

First we prove that the prefactors from Definition \ref{def:prefactors} are well-defined. 

\begin{proposition} \label{prop:En0inftyWelld}
For $n$ big enough, the prefactors $E_n^\infty$ and $E_n^0$, as defined in \eqref{eq:defEn+} and \eqref{eq:defEn-}, are well-defined, non-singular and analytic.
\end{proposition}

\begin{proof}
$E_n^0$ is obviously well-defined and analytic for any positive integer $n$. 
To see that it is non-singular, notice that any $H\in\mathcal O_F^1$ satisfies $H(n,z)^{L}=\mathcal O(n^{d-c})$ for any fixed $L$ big enough by Proposition \ref{prop:OFkGH}(a) and (e). In particular, we then have
\begin{align} \label{eq:prodhn}
(\mathbb I-H(n,z)) \left(\mathbb I+H(n,z) + H(n,z)^2 + \ldots + H(n,z)^{L-1}\right) = \mathbb I - H(n,z)^L = \mathbb I + \mathcal O(n^{d-c}),
\end{align}
uniformly for $z\in\partial D(0,n^{-a})$ as $n\to\infty$. The assumptions of Theorem \ref{lem:matching} imply that $c>d$, and thus that the right-hand side of \eqref{eq:prodhn} is invertible for $n$ big enough. Hence both factors on the left-hand side must also be invertible for $n$ big enough. In particular, $\mathbb I - H(n,z)$ is invertible for $n$ big enough. When in addition $H$ is analytic on $\overline{D(0,n^{-a})}$, the maximum modulus principle implies that $H(n,z)^L = \mathcal O(n^{d-c})$ on $\overline{D(0,n^{-a})}$. Then we can conclude that $\mathbb I - H(n,z)$ is invertible on $\overline{D(0,n^{-a})}$ for $n$ big enough.  Applying this reasoning to each of the factors $\mathbb I - (\pi^j F)^+$ in the definition \eqref{eq:defEn+} of $E_n^0$ yields that $E_n^0$ is non-singular on $\overline{D(0,n^{-a})}$ for $n$ big enough. Here we also used the assumption that $E$ is non-singular.

To prove that $E_n^\infty$ is well-defined we will use a slight abuse of notation, namely we let $(E_n^\infty)^{-1}$ be defined as the right-hand side of \eqref{eq:defEn-}, but without the inverse. To prove that $E_n^\infty$ is well-defined, it then suffices to show that $(E_n^\infty)^{-1}$ is non-singular. This would then also immediately imply that $E_n^\infty$ is non-singular. We know that each factor $\mathbb I - (\pi^j F)^-(z)$ in the definition \eqref{eq:defEn-} of $E^\infty_n(z)$ is a polynomial evaluated in $1/z$ with constant term $\mathbb I$. Then $E_n^\infty(z)^{-1}$, too, is a polynomial evaluated in $1/z$ with constant term $\mathbb I$. An application of the maximum modulus principle, to the region $1/z \in \overline{D(0,n^{-a})}$, then yields that $E_n^\infty(z)^{-1}$, as a function on $\overline{A(0;n^{-a},\infty)}$, attains its maximum on the circle $\partial D(0,n^{-a})$. The same reasoning that we used for $E_n^0$ will now prove that $(E_n^\infty)^{-1}$ is non-singular on $A(0;n^{-a},\infty)$ when $n$ is big enough. In fact, since $(E_n^\infty)^{-1}$ is non-singular for $n$ big enough, and a polynomial evaluated in $1/z$, the prefactor $E_n^\infty$ must be analytic for $n$ big enough.\\
\end{proof}

\noindent \textit{Proof of Theorem \ref{lem:matching}}. We know due to Proposition \ref{prop:0step} that
\begin{align} \label{eq:0step}
E(z) \mathring P(z) N(z)^{-1} = \mathbb I + F(z) + \mathcal O_E,
\end{align}
where $F$ is as in Definition \ref{def:F}. Trivially, $F\in\mathcal O_F^1$.  Hence Corollary \ref{cor:induction} tells us that  
\begin{align*}
\prod_{j=0}^K \left(\mathbb I - (\pi^{K-j} F)^+\right) 
E P N^{-1}
\prod_{j=0}^K \left(\mathbb I - (\pi^j F)^-\right)
= \mathbb I + \pi^{K+1} F + \mathcal O_E,
\end{align*}
with $K$ as in Definition \ref{def:prefactors}. Additionally, it tells us that $\pi^{K+1} F\in\mathcal O_F^{2^{K+1}}$. Then Proposition \ref{prop:OFkGH}(e) combined with ${2^{K+1}\geq \frac{a+c-e}{b-e}}$, tells us that $\pi^{K+1} F(z) = \mathcal O(n^{d-c})$ uniformly for $z\in\partial D(0,n^{-a})$ as $n\to\infty$. Combining this with Proposition \ref{prop:prodOAE}(a), we infer that 
\begin{align*}
\prod_{j=0}^K \left(\mathbb I - (\pi^{K-j} F)^+(z)\right) 
E(z) P(z) N(z)^{-1}
\prod_{j=0}^K \left(\mathbb I -  (\pi^{j} F)^-(z)\right)
= \mathbb I + \mathcal O(n^{d-c})
\end{align*}
uniformly for $z\in\partial D(0,n^{-a})$ as $n\to\infty$. 
Hence, we obtain the matching \eqref{eq:matchingonna} on $\partial D(0,n^{-a})$.

It remains to prove that the matching \eqref{eq:matchingonr} on $\partial D(0,r)$ is satisfied. By construction, $E_n^\infty(z)^{-1}$ equals some polynomial evaluated in $1/z$ with constant term $\mathbb I$. In particular, $E_n^\infty(z)^{-1}-\mathbb I$ is its own principal part. Then, using \eqref{eq:principalPartInt} with $\rho=n^{-a}$, we have for $z\in\partial D(0,r)$ that
\begin{align} \label{eq:EninftyCauchy}
E_n^\infty(z)^{-1} = \mathbb I - \frac{1}{2\pi i} \oint_{\partial D(0,n^{-a})} \frac{E_n^\infty(s)^{-1}-\mathbb I}{s-z} ds
\end{align}
Proposition \eqref{prop:OFkGH}(a) implies that $E_n^\infty(z)^{-1}-\mathbb I\in \mathcal O_F^1$. By Proposition \ref{prop:OFkGH}(d) any element of $\mathcal O_F^1$ is $\mathcal O\left(n^{a-b+d}\right)$ uniformly for $z\in\partial D(0,n^{-a})$ as $n\to\infty$. Applying this to \eqref{eq:EninftyCauchy} yields
\begin{align*}
E_n^\infty(z)^{-1} = \mathbb I + \mathcal O\left(n^{a-b+d} n^{-a}\right) = \mathbb I + \mathcal O\left(n^{d-b}\right)
\end{align*}
uniformly for $z\in\partial D(0,r)$ as $n\to\infty$. Since $b>d$ we may take the inverse for $n$ big enough, and we obtain \eqref{eq:matchingonr}.\qed

\begin{remark} \label{remark:globalMatching} 
I would like to point out that, while obtaining the matching in a local way seems more favorable, it is also possible to perform the matching globally. Indeed, due to the polynomial nature of $E_n^\infty(z)^{-1}$ for the constructed $E_n^\infty(z)$ in \eqref{eq:defEn-} we have $E_n^\infty(z) = \mathbb I + \mathcal O(1/z)$ as $z\to\infty$. We may now define the final transformation instead of \eqref{finalT} as
\begin{align} \nonumber
R(z) = \left\{\begin{array}{ll}
S(z) N(z)^{-1} E_n^\infty(z)^{-1}, & z\in A(0;n^{-a},\infty),\\
S(z) \mathring P(z)^{-1} E_n^0(z)^{-1}, & z\in D(0,n^{-a}).
\end{array}\right.
\end{align}
In this case we only have a matching on $\partial D(0,n^{-a})$. This approach seems to work also when there are other special points where a local parametrix has to be constructed (in case the number of special points is not infinite). An advantage is that one may even have $b\leq d$ (the assumption that $b>d$ is only used at the very end of the proof of Theorem \ref{lem:matching}). A downside may be that the expression that we could now view as the global parametrix, i.e., $E^\infty_n(z) N(z)$, would not be independent of $n$. Obtaining the matching by modifying the global parametrix has been done in \cite{DeKu, DeKuZh, KuMFWi}.\end{remark}

\pagebreak

\section{Scaling limits of correlation kernels at a special point} \label{section:scalingLimits}

\subsection{A result on scaling limits of correlations kernels}

RH analyses are often used to find scaling limits of functions of interest, e.g., functions related to orthogonal polynomials and correlation kernels. In this section we focus on scaling limits for correlation kernels. One can often (e.g. see \cite{DaKu}) express a correlation kernel $K_n$ as
\begin{align} \label{def:KY}
K_n(x,y) = \frac{1}{2\pi i(x-y)} u(y) Y(y)^{-1} Y(x) v(x),
\end{align} 
where $u(y)$ is a row vector and $v(x)$ is a column vector, and $Y$ is the ($n$-dependent) solution to a corresponding RHP. Frequently, there exists a scaling limit at some point $(x_0,y_0)\in\mathbb R^2$ of the form
\begin{align*}
\lim_{n\to\infty} \frac{1}{\mathfrak c n^b} K_n\left(x_0+\frac{x}{\mathfrak c n^b}, y_0+\frac{y}{\mathfrak c n^b}\right) = \mathbb K(x,y),
\end{align*}
where $\mathfrak c$ and $b$ are positive constants, and $\mathbb K$ is a limiting kernel that may be different for different choices of $(x_0,y_0)$. Well-known examples of limiting kernels are the sine, Airy and Bessel kernel. Frequently, one desires the scaling limit to be uniform for $x, y$ in (real) compact sets. In this section we will focus on the arguably most interesting case, namely, the case of a scaling limit at a special point. As before we assume that $0$ is a special point of our RHP, and we will investigate the scaling limit of the correlation kernel at $(x_0,y_0)=(0,0)$. In this situation, it helps to have a good understanding of
\begin{align*}
Y\left(\frac{z}{\mathfrak c n^b}\right)
\end{align*}
as $n$ becomes large. Then, as we shall see, one inevitably has to take the structure of the local parametrix around $0$ into account, which is what makes this case interesting. A succesful application of the steepest descent analysis, through transformations $Y \mapsto X \mapsto T \mapsto S \mapsto R$ (as in Section \ref{sec:setUp}), or equivalent, will create a RHP for $R$, that is trivial enough to conclude that $R$ converges to $\mathbb I$ as $n\to\infty$, 
uniformly outside the jump contour as $n\to\infty$ (e.g., see \cite{KaMcMi}). One can then invert all the transformations of the RH analysis, and say something about the large $n$ behavior of $Y$, and in particular, about potential scaling limits of $Y$. 

As mentioned in the beginning of this section, the dependence of $u(x)$ and $v(x)$ in \eqref{def:KY} on $n$ is generally not significant in the scaling regime, where $x$ and $y$ are of order $n^{-b}$. When calculating scaling limits of correlation kernels at a special point, the origin in our case, it is thus enough to understand the behavior of $Y(y_n)^{-1} Y(x_n)$, with
\begin{align} \label{eq:defxnyn}
x_n = \frac{x}{\mathfrak c n^b}, \quad\quad y_n = \frac{y}{\mathfrak c n^b}.
\end{align}
Using the same notations as in Subsection \ref{sec:setUp}, it turns out, after a considerable amount of bookkeeping, that we usually have
\begin{align*}
\frac{1}{\mathfrak c n^b} K_n(x_n,y_n) = \frac{h_n(x,y)}{2\pi i(x-y)} u_0 \Psi_+\left(n^b f(y_n)\right)^{-1} E_n^0(y_n)^{-1} R(y_n)^{-1} R(x_n) E_n^0(x_n) \Psi_+\left(n^b f(x_n)\right)  v_0,
\end{align*}
for a constant row vector $u_0$ and a constant column vector $v_0$, and some scalar factor $h_n$ that converges rapidly to $1$ as $n\to\infty$. Now we would like to argue, setting $\mathfrak c=f'(0)$ in \eqref{eq:defxnyn}, that
\begin{align} \label{eq:scalingKuv}
\lim_{n\to\infty} \frac{1}{\mathfrak c n^b} K_n\left(x_n,y_n\right)
= \frac{1}{2\pi i(x-y)} u_0 \Psi_+\left(y\right)^{-1} \Psi_+\left(x\right)  v_0,
\end{align}
uniformly for $x,y$ in compact sets. Indeed, we have $h_n(x,y)\to 1, n^b f(x_n) \to x$ and $n^b f(y_n)\to y$ as $n\to\infty$. Hence, to arrive at the limit \eqref{eq:scalingKuv}, one needs that the expression
\begin{align*}
E_n^0(y_n)^{-1} R(y_n)^{-1} R(x_n) E_n^0(x_n)
\end{align*}
is close to the unit matrix for large $n$. Indeed, the main goal of this section is to prove that such a behavior holds under reasonable assumptions. 

Let us suppose that we have indeed carried out our RH analysis, performed a double matching, and have just defined our final transformation to obtain $R$, as in \eqref{finalT}. We are then generally in the following situation. $R$ is analytic except on some oriented contour $\Sigma$, that can be partitioned as 
\begin{align} \label{eq:decompSigma}
\Sigma = \partial D(0,n^{-a}) \cup \Sigma_r \cup \partial D(0,r) \cup \Sigma_\infty.
\end{align}
$\Sigma_r$ consists of the jump curves in $A(0;n^{-a},r)$. These curves generally correspond to the lips of lenses, the corresponding jump matrices usually behave as $\mathbb I + \mathcal O(\exp(- \alpha n |z|^{\beta}))$ for some constants $\alpha, \beta > 0$ (in fact, I am confident that one has $\beta=1/b$ in natural situations). $\Sigma_\infty$ consists of the jump curves in $A(0;r,\infty)$. Theorem \ref{lem:matching} tells us that the jump matrix behaves on $\partial D(0,n^{-a})$ and $\partial D(0,r)$ as $\mathbb I + \mathcal O(n^{d-c})$ and $\mathbb I + \mathcal O(n^{b-c})$ uniformly as $n\to\infty$ respectively. Jumps in $A(0;r,\infty)$ might correspond to other local parametrices, but far away the jump matrix converges radiply to $\mathbb I$ or there are no jumps at all. If we denote the jump matrix by $\mathbb I + \Delta$, then we can summarize our situation in $\overline{D(0,r)}$ as
\begin{align} \label{eq:estimatesContours}
\Delta(z) =\left\{\begin{array}{ll}
\mathcal O(n^{d-c}), & z\in\partial D(0,n^{-a}),\\
\mathcal O(n^{b-c}), & z\in\partial D(0,r),\\
\mathcal O(\exp(- \alpha n |z|^{\beta})), & z\in \Sigma_r, 
\end{array}\right.
\end{align}
where the estimates are uniform on the indicated curves as $n\to\infty$. We then have $\lvert\lvert \Delta \rvert\rvert_{L^2(\Sigma)}\to 0$ and $\lvert\lvert \Delta \rvert\rvert_{L^\infty(\Sigma)}\to 0$ as $n\to\infty$. If, additionally, $R(z) = \mathbb I + \mathcal O(1/z)$ as $z\to\infty$, then a general theorem (e.g., see \cite{De}) implies that $R$ can be represented in the integral form
\begin{align} \label{eq:RIntegralForm}
R(z) = \mathbb I + C^\Sigma(F \Delta) = \mathbb I + \int_\Sigma \frac{F(s) \Delta(s)}{s-z} ds,
\end{align}
where $C^\Sigma$ denotes the Cauchy operator with respect to the oriented contour $\Sigma$, and $F = \mathbb I + X$, where $X$ is some function in $L^2(\Sigma)$ that satisfies $\lvert\lvert X \rvert\rvert_{L^2(\Sigma)} \to 0$ as $n\to\infty$. There is frequently an estimate like
\begin{align} \label{eq:estimateRtoI}
R(z) = \mathbb I + \mathcal O\left(\frac{1}{n^\delta (1+|z|)}\right)
\end{align}
uniformly for $z\in\mathbb C\setminus \Sigma$ as $n\to\infty$, for some fixed $\delta>0$ (e.g., see Lemma B.0.2. in \cite{KaMcMi} or the proof of Proposition 2.5.1 in \cite{BlLi}). A naive application of \eqref{eq:estimateRtoI}, using a standard argument with Cauchy's integral formula, combined with results on the behavior of $E_n^0$ (see Lemma \ref{thm:scalingb}), then yields that 
\begin{align} \label{eq:naive}
E_n^0(y_n)^{-1} R(y_n)^{-1} R(x_n) E_n^0(x_n) = \mathbb I + \mathcal O(n^{e-b}(x-y)) + \mathcal O(n^{a-b+d-\delta}(x-y))
\end{align}
uniformly for $x,y$ in compact sets as $n\to\infty$. It would be nice if $\delta > a-b+d$, but there generally is no a priori reason that this is the case. Therefore we include Theorem \ref{thm:corKerScaling} below. 

Actually, our situation is trickier than in \cite{De}, \cite{KaMcMi} and \cite{BlLi}, because our contour $\Sigma$ varies with $n$. However, under some extra conditions we may derive similar properties for $R$ by arguments along the lines of Appendix A from \cite{BlKu}. We shall not worry about such subtleties  and just take \eqref{eq:estimateRtoI}, or rather the weaker statement that $R\to\mathbb I$ uniformly as $n\to\infty$, and \eqref{eq:RIntegralForm} as a starting point. 

\begin{theorem} \label{thm:corKerScaling}
Suppose that the conditions of Theorem \ref{lem:matching} are met and let $E_n^0$ be as in Definition \ref{def:prefactors}. Assume, additionally, that $C$ is  uniformly bounded on $\partial D(0,n^{-e})$ and that 
\begin{align} \label{eq:condOncade}
c\geq \min\left(\frac{3}{2}a+d,\frac{3}{2}a + 2d - e\right).
\end{align}
Suppose that $R:\mathbb C\setminus\Sigma\to\mathbb C^{m\times m}$ converges uniformly to $\mathbb I$ as $n\to\infty$, and that it can be written in the integral form \eqref{eq:RIntegralForm}, where 
\begin{itemize}
\item[(i)] $\Sigma$ is an oriented contour as in \eqref{eq:decompSigma}, where $\Sigma_r$ and $\Sigma_\infty$ consist of a finite union of smooth curves, and the inversion $s\mapsto s^{-1}$ is bounded in $L^2(\Sigma_\infty)$ sense by some ($n$-independent) constant. 
\item[(ii)] $\Delta\in L^2(\Sigma)\cap L^\infty(\Sigma)$ satisfies the estimates \eqref{eq:estimatesContours} for some $\alpha>0$ and some $0<\beta<1/a$, and both $\lvert\lvert \Delta \rvert\rvert_{L^2(\Sigma)}\to 0$ and $\lvert\lvert \Delta \rvert\rvert_{L^\infty(\Sigma)}\to 0$ as $n\to\infty$. 
\item[(iii)] $F=\mathbb I + X$, with $X\in L^2(\Sigma)$ and $\lvert\lvert X \rvert\rvert_{L^2(\Sigma)}\to 0$ as $n\to\infty$.
\end{itemize}
Then we have, with $x_n$ and $y_n$ as in \eqref{eq:defxnyn}, that
\begin{align} \label{eq:ERREerror}
E_n^0(y_n)^{-1} R(y_n)^{-1} R(x_n) E_n^0(x_n) = \mathbb I + \mathcal O\left(n^{d-b}(x-y)\right)  + \mathcal O\left(n^{e-b}(x-y)\right) 
\end{align}
uniformly for $x,y$ in compact sets as $n\to\infty$.
\end{theorem}
In the natural situation \eqref{eq:naturalC}, that is,
\begin{align*}
C(z) = \frac{z}{f(z)}\left(C_1  +\frac{C_2}{(n^b z)} \frac{z}{f(z)} + \ldots \frac{C_k}{(n^b z)^{k-1}} \frac{z^{k-1}}{f(z)^{k-1}} \right),
\end{align*}
we see that $C$ is uniformly bounded on all of $\overline{A(0;n^{-b},r)}$ as $n\to\infty$. Assuming that $b>e$, we then certainly have that $C$ is uniformly bounded on $\partial D(0, n^{-e})$ as $n\to\infty$. 

\begin{remark} \label{remark:cBigEnough} 
I would like to emphasize an important insight, which is that the condition \eqref{eq:condOncade} can always be arranged in the natural situation \eqref{eq:naturalC}, namely, by taking enough terms in the asymptotic expansion \eqref{eq:behavPsi}. Indeed, we have that $c=(b-a)(k+1)$ in general for any fixed $k$.\\
If necessary, e.g., in cases were $C$ might have a more complicated description, one can omit the condition on $c$ altogether, but then \eqref{eq:ERREerror} has to be replaced by
\begin{align*}
\hspace{-0.5cm}E_n^0(y_n)^{-1} R(y_n)^{-1} R(x_n) E_n^0(x_n) = \mathbb I + \mathcal O(n^{d-b}(x-y)) + \mathcal O(n^{e-b}(x-y)) 
+ \mathcal O(n^{\frac{3}{2}a-b-c+2d}(x-y)).
\end{align*}
A meaningful application of Theorem \ref{thm:corKerScaling} would then nevertheless require the condition $c> \frac{3}{2}a-b+2d$. 
\end{remark}

\subsection{Behavior of the prefactor close to the origin}

In order to prove Theorem \ref{thm:corKerScaling}, one thing we need to understand is how the prefactor $E_n^0$ behaves close to the origin. The following lemma yields estimates for the analytic prefactor near the origin. 

I would like to point out that this result may also be useful for the calculation of scaling limits that do \textit{not} concern correlation kernels. 

\begin{lemma} \label{thm:scalingb}
Suppose that the assumptions of Theorem \ref{lem:matching} are met and that, additionally, $C$ is uniformly bounded on $\partial D(0,n^{-e})$. Then the prefactor $E_n^0$, as in Definition \ref{def:prefactors}, satisfies uniformly for $z\in\overline{D(0,n^{-e})}$ that as $n\to\infty$
\begin{align} \label{eq:E+E0inb}
E_n^0(z) = E(0) \left(\mathbb I +\mathcal O(n^{e-b}) + \mathcal O\left(n^{e} z\right)\right)
\hspace{0.2cm}\text{and}\quad 
E_n^0(z)^{-1} =  \left(\mathbb I +\mathcal O(n^{e-b}) + \mathcal O\left(n^{e} z\right)\right) E(0)^{-1}.
\end{align}
Furthermore, for any fixed $0<\rho<1$, we have uniformly for $z,w\in \overline{D(0, n^{-e} \rho)}$ that as $n\to\infty$
\begin{align} \label{eq:EzEwInez-w}
E_n^0(z)^{-1} E_n^0(w) = \mathbb I + \mathcal O(n^e (z-w)).
\end{align}
\end{lemma}

\begin{proof} By the assumption \eqref{eq:assumpE} of Theorem \ref{lem:matching} we have uniformly for $z,w\in \partial D(0,n^{-a})$ that
\begin{align} \label{eq:EEe}
E(z)^{-1} E(w) = \mathbb I + \mathcal O(n^e(z-w))
\end{align}
as $n\to\infty$. We claim that \eqref{eq:EEe} is actually true uniformly for $z,w\in\overline{D(0,n^{-a})}$ as $n\to\infty$. To prove this, we write \eqref{eq:EEe} more suggestively as
\begin{align}  \label{eq:EEe2}
n^{-e} E(z)^{-1} \frac{E(z)-E(w)}{z-w}  = \mathcal O(1)
\end{align}
uniformly for $z,w\in \partial D(0,n^{-a})$ as $n\to\infty$. For fixed $n$, the function on the left-hand side is analytic in both of its variables (the singularity in $z=w$ is removable). All its $m^2$ entries are, in absolute value, bounded by some ($n$-independent) positive constant, implied by the $\mathcal O(1)$ term. We can now apply the maximum modulus principle, with respect to both variables, to each of the entries, for fixed $n$. Since $n$ is arbitrary in this argument, we infer that \eqref{eq:EEe2} is actually true uniformly on $\overline{D(0,n^{-a})}$ as $n\to\infty$. Rewriting \eqref{eq:EEe2} then yields the claim, i.e.,
\begin{align} \label{eq:EEeClaim}
E(z)^{-1} E(w) &= \mathbb I + \mathcal O(n^e(z-w)), & \text{uniformly for }z,w\in\overline{D(0,n^{-a})}\text{ as }n\to\infty.
\end{align}

It follows from the claim \eqref{eq:EEeClaim} and $e\geq a$ that uniformly for $z,w\in \overline{D(0, n^{-e})}$ as $n\to\infty$
\begin{align} \label{eq:EEee}
E(z)^{-1} E(w) 
= \mathbb I + \mathcal O(n^e (z-w))
= \mathcal O(1).
\end{align}
Let $l$ be any positive integer. That $C$ is uniformly bounded on $\partial D(0,n^{-e})$, combined with \eqref{eq:EEee} and Definition \ref{def:F}, implies that uniformly for $z_1,\ldots,z_l\in\partial D(0,n^{-e})$ as $n\to\infty$
\begin{align} \nonumber
F(z_1) \cdots F(z_l) &= (n^b z)^{-l} E(z_1) C(z_1)\\ \nonumber
& \hspace{2cm} \left(E(z_1)^{-1} E(z_2)\right) C(z_2)  \cdots 
C(z_{l-1}) \left(E(z_{l-1})^{-1} E(z_l)\right) C(z_l) E(z_l)^{-1}\\ \nonumber
&= E(0) \left(E(0)^{-1} E(z_1)\right) \mathcal O\left(n^{-(b-e) l}\right) \left(E(z_l)^{-1} E(0)\right) E(0)^{-1}\\ \label{eq:proofAAe}
&= E(0) \mathcal O\left(n^{-(b-e)}\right) E(0)^{-1}.
\end{align}
We used the assumption $b>e$ in the last step. 
We conclude that any term of the form \eqref{eq:piFsumint3} as in Definition \ref{def:Ln} equals $E(0) \mathcal O\left(n^{-(b-e)}\right) E(0)^{-1}$ uniformly on $\partial D(0,n^{-e})$, which follows by taking the integrations over $\partial D(0,n^{-e})$. We also have that $E_n^0 E^{-1} - \mathbb I \in\mathcal O_F^1$, which follows from \eqref{eq:defEn+}, \text{Definition \ref{def:OFk}} and Proposition \ref{prop:OFkGH}(a) and (b). Combining these two facts yields that
\begin{align} \nonumber
E_n^0(z) &= \left(E_n^0(z) E(z)^{-1} - \mathbb I\right) E(z) + E(z)\\ \nonumber
&= E(0) \mathcal O\left(n^{-(b-e)}\right) E(0)^{-1} E(z) + E(z)\\ \label{eq:forInverse}
&= E(0) \left(\mathcal O\left(n^{-(b-e)}\right) + \mathbb I\right) E(0)^{-1} E(z)
\end{align}
uniformly for $z\in \partial D(0,n^{-e})$ as $n\to\infty$. The maximum modulus principle, applied to $E(0)^{-1} E_n^0(z) E(z)^{-1} E(0)-\mathbb I$ then yields that \eqref{eq:forInverse} is actually true uniformly for $z\in\overline{D(0,n^{-e})}$ as $n\to\infty$. Combining this with \eqref{eq:EEee} (for $w=0$) then gives that
\begin{align}\label{eq:En+ongel}
E_n^0(z) = E(0) \left(\mathbb I +\mathcal O\left(n^{e-b}\right) + \mathcal O(n^{e} z)\right)
\end{align}
 uniformly for $z\in\overline{D(0,n^{-e})}$ as $n\to\infty$. To prove the second estimate in \eqref{eq:E+E0inb}, we use \eqref{eq:forInverse}, to see that
\begin{align} \nonumber
E_n^0(z)^{-1} &= E(z)^{-1} E(0) \left(\mathbb I + \mathcal O(n^{e-b})\right)^{-1} E(0)^{-1}\\ \nonumber
&= \left(\mathbb I + \mathcal O(n^e z)\right) \left(\mathbb I + \mathcal O(n^{e-b})\right) E(0)^{-1}\\ \label{eq:forInverse2}
&= \left(\mathbb I +\mathcal O(n^{e-b}) + \mathcal O(n^{e} z)\right) E(0)^{-1}
\end{align}
uniformly for $z\in\overline{D(0,n^{-e})}$ as $n\to\infty$. Here we have used the claim \eqref{eq:EEee} in the second line (for $w=0$). 

Let us now prove the second part of the theorem. Let $0<\rho<1$ and suppose that $z,w\in \overline{D(0,n^{-e} \rho)}$. By Cauchy's integral formula and \eqref{eq:En+ongel} we have that
\begin{align} \nonumber
E_n^0(w) - E_n^0(z) &= \frac{z-w}{2\pi i} \oint_{\partial D(0,n^{-e})} \frac{E_n^0(s)-E_n^0(0)}{(s-z)(s-w)} ds\\ \nonumber
&= E(0) \frac{z-w}{2\pi i} \oint_{\partial D(0,n^{-e})} \frac{\mathcal O(n^{e-b})+\mathcal O(n^{e} s)}{(s-z)(s-w)} ds\\ \label{eq:EEEEE}
&= E(0) \frac{|z-w|}{(1-\rho)^2} \mathcal O(n^e)
= E(0) \mathcal O(n^e(z-w))
\end{align}
uniformly for $z,w\in \overline{D(0,n^{-e}\rho)}$ as $n\to\infty$. Here we have used that $e<b$. Combining \eqref{eq:EEEEE} with \eqref{eq:forInverse2}, we have uniformly for $z,w\in \overline{D(0, n^{-e}\rho)}$ that as $n\to\infty$
\begin{align} \nonumber
E_n^0(z)^{-1} E_n^0(w) &= \mathbb I + E_n^0(z)^{-1} (E_n^0(w)-E_n^0(z))\\ \nonumber
&= \mathbb I + \left(\mathbb I +\mathcal O(n^{e-b}) + \mathcal O(n^{e} z)\right) E(0)^{-1} E(0) \mathcal O(n^e (w-z))\\ \nonumber
&= \mathbb I + \mathcal O(n^e (z-w)).
\end{align}
We used here that $\mathbb I +\mathcal O(n^{e-b}) + \mathcal O(n^{e} z)=\mathcal O(1)$ uniformly for $z,w\in \overline{D(0, n^{-e} \rho)}$ as $n\to\infty$.
\end{proof}

The condition that \eqref{eq:EzEwInez-w} is true provided that $0<\rho<1$, is more or less irrelevant in the natural situation \eqref{eq:naturalC}, where we have the freedom to move our estimates to slightly bigger circles. We may in such cases essentially set $\rho$ equal to $1$, or any other positive number for that matter. 

One way to interpret Lemma \ref{thm:scalingb} is that $E_n^0$ takes over the estimates \eqref{eq:assumpE} of $E$ in Theorem \ref{lem:matching} on $\overline{D(0,n^{-e})}$. The theorem should be convenient when one is interested in the behavior of $E_n^0(z)$ in the scaling regime, i.e., for $z$ with modulus of order $n^{-b}$.

Lemma \ref{thm:scalingb} does not provide any information on the behavior of $E(0)$. We mention that $E(0)$ is, in general, a function of $n$. I suspect that there is no general result possible concerning the behavior of $E(0)$, its properties will depend strongly on the specifics of the particular RHP. In situations where one really needs the behavior of $E(0)$, for scaling limits of orthogonal polynomials perhaps, one should therefore do some extra work. The value of $E(0)$ will not be an issue when we are studying correlation kernels however, as we shall see. 

\subsection{Proof of Theorem \ref{thm:corKerScaling}}

To prove Theorem \ref{thm:corKerScaling} we also need to be precise about the behavior of $R$ near the origin. There is generally an estimate like \eqref{eq:estimateRtoI} for some fixed $\delta>0$ as $n\to\infty$ uniformly, but this might not be sharp enough for our purposes. For this reason we include the following lemma. 

\begin{lemma} \label{prop:R-R}
Under the asumptions of Theorem \ref{thm:corKerScaling} we have that
\begin{align}
R(y_n)^{-1} R(x_n) = \mathbb I + \mathcal O(n^{- b}(x-y)) + \mathcal O(n^{\frac{3}{2}a - b - c + d}(x-y)),
\end{align}
uniformly for $x,y$ in compact sets that as $n\to\infty$.
\end{lemma}

\begin{proof}
So let us consider $x,y$ in a compact set. As before, we denote by $C^\sigma$ the Cauchy-operator with respect to a contour $\sigma$. We notice, using the first estimate in \eqref{eq:estimatesContours} and Cauchy-Schwarz, that as $n\to\infty$
\begin{align} \nonumber
& C^{\partial D(0,n^{-a})}(F \Delta)(x_n) -  C^{\partial D(0,n^{-a})}(F \Delta)(y_n)
= \frac{x_n-y_n}{2\pi i} \oint_{|s|=n^{-a}} \frac{F(s) \Delta(s)}{(s-x_n)(s-y_n)} ds\\ \nonumber
&= \mathcal O\left(n^{-b}(x-y) \quad n^{2a} \quad \lvert\lvert F \rvert\rvert_{L^2(\partial D(0,n^{-a}))} 
\quad \lvert\lvert \Delta\rvert\rvert_{L^2(\partial D(0,n^{-a}))}\right)\\ \nonumber
&= \mathcal O\left(n^{2a-b}(x-y) \quad (2\pi n^{-a} m+\lvert\lvert X \rvert\rvert_{L^2(\Sigma)}) \quad \lvert\lvert \Delta\rvert\rvert_{L^\infty(\partial D(0,n^{-a})} n^{-a/2} \right)\\ \label{int1}
&= \mathcal O(n^{\frac{3}{2}a-b-c+d}(x-y)).
\end{align}
Here we have used that $2\pi n^{-a} m+\lvert\lvert X \rvert\rvert_{L^2(\Sigma)}$ is uniformly bounded as $n\to\infty$, which is a consequence of the fact that $\lvert\lvert X \rvert\rvert_{L^2(\Sigma)}\to 0$ as $n\to\infty$. Analogously, using the second estimate in \eqref{eq:estimatesContours} and Cauchy-Schwarz, we have
\begin{align} \nonumber
C^{\partial D(0,r)}(F \Delta)(x_n) &-  C^{\partial D(0,r)}(F \Delta)(y_n)
= \frac{x_n-y_n}{2\pi i} \oint_{|s|=r} \frac{F(s) \Delta(s)}{(s-x_n)(s-y_n)} ds\\ \label{int2}
&= \mathcal O(n^{-b}(x-y) \quad n^{d-b}) = \mathcal O(n^{- b}(x-y)),
\end{align}
and, using the third estimate in \eqref{eq:estimatesContours} and Cauchy-Schwarz,
\begin{align} \nonumber
C^{\Sigma_r}(F \Delta)(x_n) &-  C^{\Sigma_r}(F \Delta)(y_n)
= \frac{x_n-y_n}{2\pi i} \int_{\Sigma_r} \frac{F(s) \Delta(s)}{(s-x_n)(s-y_n)} ds\\ \label{int3}
&= \mathcal O\left(n^{-b}(x-y) \quad n^{2 a} \quad \exp(-\alpha n^{1-a \beta})\right) = \mathcal O(n^{- b}(x-y)),
\end{align}
as $n\to\infty$. Here we also used the assumptions $\alpha>0$ and $\beta<1/a$. Lastly, we have
\begin{align} \nonumber
C^{\Sigma_\infty}(F \Delta)(x_n) &-  C^{\Sigma_\infty}(F \Delta)(y_n)
= \frac{x_n-y_n}{2\pi i} \int_{\Sigma_\infty} \frac{F(s) \Delta(s)}{(s-x_n)(s-y_n)} ds\\ \label{eq:inftycontour}
&= \frac{x_n-y_n}{2\pi i} \int_{\Sigma_\infty} \frac{\Delta(s)}{(s-x_n)(s-y_n)} ds + \frac{x_n-y_n}{2\pi i} \int_{\Sigma_\infty} \frac{X(s) \Delta(s)}{(s-x_n)(s-y_n)} ds
\end{align}
By the assumption that the inversion $s\mapsto s^{-1}$ is uniformly bounded in $L^2(\Sigma_\infty)$ sense for all $n$, and that $\lvert\lvert\Delta\rvert\rvert_{L^\infty(\Sigma)}\to 0$ as $n\to\infty$, we have that 
\begin{align} \label{eq:estContourInfty1}
\frac{x_n-y_n}{2\pi i} \int_{\Sigma_\infty} \frac{\Delta(s)}{(s-x_n)(s-y_n)} ds
= \mathcal O(n^{-b}(x-y))
\end{align}
uniformly as $n\to\infty$. The boundedness of $s\mapsto (s-x_n)^{-1} (s-y_n)^{-1}$ outside $D(0,r)$, combined with the assumptions $\lvert\lvert X\rvert\rvert_{L^2(\Sigma)}\to 0$ and $\lvert\lvert\Delta\rvert\rvert_{L^2(\Sigma)}\to 0$ as $n\to\infty$ and Cauchy-Schwarz, yields that
\begin{align} \label{eq:estContourInfty2}
\frac{x_n-y_n}{2\pi i} \int_{\Sigma_\infty} \frac{X(s) \Delta(s)}{(s-x_n)(s-y_n)} ds = \mathcal O(n^{-b}(x-y))
\end{align}
uniformly as $n\to\infty$. Then \eqref{eq:estContourInfty1} and \eqref{eq:estContourInfty2} imply that
\begin{align} \label{int4}
C^{\Sigma_\infty}(F \Delta)(x_n) &-  C^{\Sigma_\infty}(F \Delta)(y_n) = \mathcal O(n^{-b}(x-y))
\end{align}
uniformly as $n\to\infty$. Putting it all together, i.e., using \eqref{int1}, \eqref{int2}, \eqref{int3} and \eqref{int4}, we must conclude that 
\begin{align} \nonumber
R(x_n) - R(y_n) = \mathcal O(n^{\frac{3}{2}a-b-c+d}(x-y)) + \mathcal O(n^{-b}(x-y))
\end{align}
uniformly as $n\to\infty$. Hence uniformly as $n\to\infty$
\begin{align} \nonumber
R(y_n)^{-1} R(x_n) = \mathbb I + R(y_n)^{-1} (R(x_n) - R(y_n))
= \mathbb I + \mathcal O(n^{\frac{3}{2}a-b-c+d}(x-y)) + \mathcal O(n^{- b}(x-y)).
\end{align}
Here we have used the assumption that $R\to\mathbb I$ uniformly as $n\to\infty$, which implies that $R^{-1}$ is uniformly bounded as $n\to\infty$.\\
\end{proof}

We remark that the assumption on $E_n^0$ and the asumption \eqref{eq:condOncade} on $c$ in Theorem \ref{thm:corKerScaling} did not play a role in the proof of Lemma \ref{prop:R-R}. We are now ready to give the proof for Theorem \ref{thm:corKerScaling}.\\

\noindent\textit{Proof of Theorem \ref{thm:corKerScaling}}. By Lemma \ref{prop:R-R} we have that
\begin{align} \nonumber
E_n^0(y_n)^{-1} R(y_n)^{-1} R(x_n) E_n^0(x_n)
&= E_n^0(y_n)^{-1} E_n^0(x_n)
+ E_n^0(y_n)^{-1} \mathcal O(n^{-b}(x-y)) E_n^0(x_n)\\ \label{eq:ERRE}
&\quad + E_n^0(y_n)^{-1} \mathcal O(n^{\frac{3}{2}a-b-c+d}(x-y)) E_n^0(x_n).
\end{align}
uniformly for $x,y$ in compact sets as $n\to\infty$. Due to the assumptions of Theorem \ref{lem:matching} that $E$ and $E^{-1}$ are $\mathcal O(n^\frac{d}{2})$ uniformly on $\partial D(0,n^{-a})$ as $n\to\infty$, and the maximum modulus principle, we have that $E(0)$ and $E(0)^{-1}$ are $\mathcal O(n^\frac{d}{2})$ as $n\to\infty$. Combining this with the first estimates \eqref{eq:E+E0inb} of Lemma \ref{thm:scalingb} yields 
\begin{align*}
E_n^0(y_n)^{-1} = \mathcal O(n^\frac{d}{2})\quad\quad\text{and}\quad\quad
E_n^0(x_n) = \mathcal O(n^\frac{d}{2})
\end{align*}
uniformly for $x,y$ in compact sets as $n\to\infty$. Hence \eqref{eq:ERRE} turns into 
\begin{align*}
E_n^0(y_n)^{-1} R(y_n)^{-1} R(x_n) E_n^0(x_n)
= E_n^0(y_n)^{-1} E_n^0(x_n) + \mathcal O(n^{d-b}(x-y)) + \mathcal O(n^{\frac{3}{2}a-b-c+2d}(x-y))
\end{align*}
uniformly for $x,y$ in compact sets as $n\to\infty$. Now invoking the second estimate \eqref{eq:EzEwInez-w} of Lemma \ref{thm:scalingb}, we arrive at
\begin{align*}
E_n^0(y_n)^{-1} R(y_n)^{-1} R(x_n) E_n^0(x_n)
&= \mathbb I + \mathcal O(n^{e-b}(x-y)) + \mathcal O(n^{d-b}(x-y)) + \mathcal O(n^{\frac{3}{2}a-b-c+2d}(x-y))\\
&= \mathbb I + \mathcal O(n^{d-b}(x-y)) + \mathcal O(n^{e-b}(x-y))
\end{align*}
uniformly for $x,y$ in compact sets as $n\to\infty$. We used the assumption $c\geq \min(\frac{3}{2}a+d,\frac{3}{2}a+2d-e)$ in the last step.
\qed

\pagebreak

\section{Examples} \label{sec:examples}

In this section we apply the main result, Theorem \ref{lem:matching}, to some examples from the literature. I hope that these clarify why we expect the situation of \eqref{eq:matchingAsympS}, and why we expect that $\mathring P$ is of the form \eqref{eq:defmathringP}. We repeat the last equation for convenience.
\begin{align*}
\mathring P(z) = \mathring E(z) \Psi\left(n^b f(z)\right) D(z) e^{n D_\varphi(z)}.
\end{align*}
The first example that we treat is worked out in detail. The example is taken from \cite{KuMo}, the article that was the main inspiration for the present one. It should serve as a guiding example to those who intend to apply the main result. In the example, we use the notations of the current paper, rather than the notations of \cite{KuMo}. 

In the remaining two examples we will be more brief. It is not our intention to be fully rigorous here, but rather to show the strength of the main result. In these two examples, we shall identify $\mathring E, \Psi, f, D, D_\varphi$, as well as $a,b,c,d,e,$ with their counterparts from the corresponding articles. The interested reader will have to consult these articles in order to fully comprehend every detail. 

\subsection{Muttalib-Borodin ensemble with parameter $\theta=\frac{1}{2}$}

The article \cite{KuMo} about the local universality at the hard edge of the Muttalib-Borodin ensemble with $\theta=\frac{1}{2}$ has been the main inspiration for this work. The goal was to find a certain scaling limit of a correlation kernel $K_n(x,y)$ at the origin. The paper uses a $3\times 3$ RHP that depends on a parameter $\alpha>-1$ and some external field $V(x)$ that, among other requirements, is analytic around $x=0$. It is the first article where a double matching was used, although the iteration steps of the current paper are more efficient than those used in \cite{KuMo}, as we shall see. We will not be precise about the specifics of the ensemble, or the relation between the correlation kernel and the initial RHP. Instead, we focus on the corresponding matching problem. 

\begin{figure}[t]
\begin{center}
%
%
\resizebox{10.5cm}{6cm}{%
\begin{tikzpicture}[>=latex]
	\draw[-] (-7,0)--(7,0);
	\draw[-] (-3,-4)--(-3,4);
	\draw[fill] (-3,0) circle (0.1cm);
	\draw[fill] (3,0) circle (0.1cm);
	
	\node[above] at (-3.25,0) {\large 0};	
	\node[above] at (3.1,0.05) {\large $q$};	
	\node[above] at (-3.5,1.8) {\large $\Delta_1^+$};	
	\node[below] at (-3.5,-1.8) {\large $\Delta_1^-$};	
	\node[above] at (-5,0.05) {\large $\Delta_2$};	
	\node[above] at (0,0.05) {\large $\Delta_1$};
	
	\draw[->, ultra thick] (-3,0) -- (-3,1.2);
	\draw[-, ultra thick] (-3,0) -- (-3,2);
	
	\draw[->, ultra thick] (-3,0) -- (-3,-1.2);
	\draw[-, ultra thick] (-3,0) -- (-3,-2);
	
	\draw[->, ultra thick] (3,0) -- (5,0);
	\draw[->, ultra thick] (-3,0) -- (0,0);
	\draw[-, ultra thick] (-3,0) -- (7,0);
	
	\draw[->, ultra thick] (-7,0) -- (-4.5,0);
	\draw[-, ultra thick] (-7,0) -- (-3,0);
	
	\draw[-, ultra thick] (-3,2) to [out=0, in=135] (3,0);
	\draw[-, ultra thick] (-3,-2) to [out=0, in=225] (3,0);

	\draw[->, ultra thick] (-0.9,1.88) to (-0.7,1.85);
	\draw[->, ultra thick] (-0.9,-1.88) to (-0.7,-1.85);
	
\end{tikzpicture}
}
%
\caption{Contour $\Sigma_{S} = \mathbb{R} \cup \Delta_{1}^{\pm}$ for the RHP for $S$ (taken from \cite{KuMo}).\label{FigS}}
\end{center}
\end{figure}
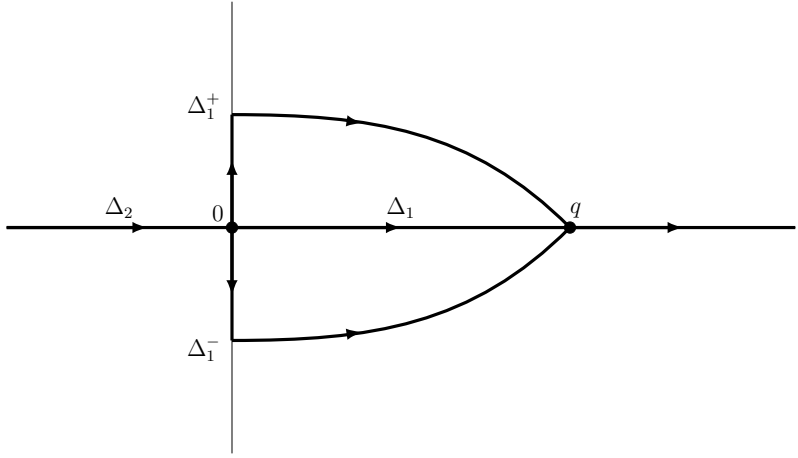

In what follows we will take $\omega=e^\frac{2\pi i}{3}$ and $\beta=\alpha+\frac{1}{4}$. After performing the standard transformations in the steepest descent analysis, we have a RHP for a function $S$, which has jumps on a contour $\Sigma_S$ that consists of five curves $\Delta_1, \Delta_2, (q,\infty)$ and $\Delta_1^\pm$, as in Figure \ref{FigS}. The curves $\Delta_1^\pm$ correspond to the lips of a lens. For $\alpha>-\frac{1}{2}$ the RHP satisfied by $S$ takes the following form.
\begin{itemize}
\item[RH-S1] $S : \mathbb{C}\setminus\Sigma_{S} \to \mathbb{C}^{3\times 3}$ is analytic. 
\item[RH-S2] On $\Sigma_S$ we have the jumps
\begin{align*} 
S_{+}(x) &= S_{-}(x)\begin{pmatrix} 0 & x^{\beta} & 0\\ -x^{-\beta} & 0 & 0\\ 0 & 0 & 1\end{pmatrix}, && x\in\Delta_{1}, \\
S_{+}(z) &= S_{-}(z) \begin{pmatrix} 1 & 0 & 0\\ z^{-\beta} e^{2 n\varphi_1(z)} & 1 & 0\\ 0 & 0 & 1\end{pmatrix}, && z\in\Delta_{1}^{+} \cup \Delta_1^-,\\ 
S_{+}(x) &= S_{-}(x) \begin{pmatrix} 1 & x^{\beta} e^{-2 n\varphi_1(x)} & 0\\ 0 & 1 & 0\\ 0 & 0 & 1\end{pmatrix}, && x\in (q,\infty), \\
S_{+}(x) &= S_{-}(x) \begin{pmatrix} 1 & 0 & 0\\ 0 & 0 & -1\\ 0 & 1 & 0\end{pmatrix}, && x\in\Delta_{2}, 
\end{align*}
\item[RH-S3] As $z \to\infty$
\begin{align*}
S(z) &= \left(\mathbb{I}+\mathcal{O}\left(\frac{1}{z}\right)\right)
\begin{pmatrix} 1 & 0 & 0\\ 0 & z^{\frac{1}{4}} & 0\\ 0 & 0 & z^{-\frac{1}{4}}\end{pmatrix}
\begin{pmatrix} 1 & 0 & 0\\ 0 & \frac{1}{\sqrt{2}} & \frac{i}{\sqrt{2}}\\ 0 & \frac{i}{\sqrt{2}} & \frac{1}{\sqrt{2}}\end{pmatrix}.
\end{align*}
\item[RH-S4] As $z\to 0$
\begin{align*}
S(z) &=  \mathcal{O}\begin{pmatrix} z^{-\alpha-\frac{1}{2}}   & z^{-\frac{1}{4}}   & z^{-\frac{1}{4}}  \\ z^{-\alpha-\frac{1}{2}}   & z^{-\frac{1}{4}}   & z^{-\frac{1}{4}}  \\ z^{-\alpha-\frac{1}{2}}   & z^{-\frac{1}{4}}   & z^{-\frac{1}{4}}  \end{pmatrix}
&& \text{ for } z \text{ inside the lense,}\\
S(z) &=  \mathcal{O}\begin{pmatrix} 1 & z^{-\frac{1}{4}}   & z^{-\frac{1}{4}}  \\ 1 & z^{-\frac{1}{4}}   & z^{-\frac{1}{4}}  \\ 1 & z^{-\frac{1}{4}}   & z^{-\frac{1}{4}}  \end{pmatrix}
&& \text{ for } z \text{ outside the lense.}
\end{align*}
\end{itemize}

Here $\varphi_1$ in RH-S2 is a $\varphi$-function having to do with the normalization and opening of the lens (see (5.13) in \cite{KuMo}). For $-1<\alpha\leq -\frac{1}{2}$, the asymptotic behavior in RH-S4 takes a different form. The particular form of RH-S4, be it for $\alpha>-\frac{1}{2}$ or be it for $-1<\alpha\leq \frac{1}{2}$, will be irrelevant however, in what follows. 

We have a global parametrix $N$ (see Section 4.4 in \cite{KuMo}) that only has jumps on $\Delta_1$ and $\Delta_2$, the same jumps as in RH-S2. Additionally, $N$ satisfies the same asymptotics as in RH-S3. The specifics of $N$ will not be relevant however. 
The local parametrix problem around $z=0$ takes the following form.
\begin{itemize}
\item[RH-P1] $P$ is analytic on $\overline{D(0,r)}\setminus \Sigma_S$.
\item[RH-P2] $P$ has the same jumps as $S$ has on $\Sigma_S \cap \overline{D(0,r)}$. 
\item[RH-P3] $P$ has the same behavior as $S$ has as $z\to 0$. 
\end{itemize}
Usually, one adds the matching condition as RH-P4. Since we plan to get a double matching instead of an ordinary matching, we leave it out. The initial construction for the local parametrix has the form 
\begin{align} \label{mutBorP}
\mathring P(z) = \mathring E_n(z) \Psi\left(n^3 f(z)\right) D(z) e^{n D_\varphi(z)},
\end{align}
where the expressions are defined as follows. 
\begin{align*}
D(z) = \begin{pmatrix} 1 & 0 & 0\\ 0 & z^\beta & 0\\ 0 & 0 & z^\beta\end{pmatrix},
\end{align*}
\begin{align*}
D_\varphi(z) = \frac{2}{3} \begin{pmatrix} 
2 \varphi_1(z)+\varphi_2(z) & 0 & 0\\ 
0 & \varphi_2(z)-\varphi_1(z) & 0
\\ 0 &  0 & \varphi_1(z)+2\varphi_2(z)
\end{pmatrix}.
\end{align*}
Here $\varphi_2$ is also a $\varphi$-function that originated from the normalization and opening the lens (see (5.14) in \cite{KuMo}).  In general, the number of relevant $\varphi$-functions equals the number of $g$-functions that were used in the normalization. The conformal map $f$ is constructed using the $\varphi$-functions, namely
\begin{align*}
f(z) = \frac{8}{729} \times \left\{\begin{array}{rl} \left(\omega^2 \varphi_1(z) - \varphi_2(z)\right)^3, & \operatorname{Im}(z)>0,\\
\left(\omega \varphi_1(z) - \varphi_2(z)\right)^3, & \operatorname{Im}(z)<0. \end{array}\right.
\end{align*}
It is proved in Proposition 5.6 in \cite{KuMo} that $f$ is indeed a conformal map that maps $0$ to itself, and that maps positive numbers to positive numbers. Then $f$ respects the jump contour of $S$, i.e., $f(\Delta_1\cap D(0,r))\subset \mathbb R^+, f(\Delta_1^+\cap D(0,r))\subset i\mathbb R^+$, etcetera, and the orientations are preserved. 

To understand $\Psi$ we have to introduce a special function. Namely, we use the Meijer G-function
\begin{align*}
G_{3,0}^{0,3}\left(\left. \begin{array}{c} -\\ 0,-\alpha, - \alpha - \frac{1}{2}\end{array} \right| z\right)
= \frac{1}{2\pi i} \int_{L} \Gamma(s) \Gamma(s-\alpha) \Gamma(s-\alpha - \tfrac{1}{2}) z^{-s} ds,
\end{align*}
where $L$ encircles the interval $(-\infty, \max(0,\alpha+\tfrac{1}{2})]$. It is a solution to the linear differential equation 
\begin{align*}
\vartheta(\vartheta+\alpha)(\vartheta+\alpha+\frac{1}{2}) \phi + z \phi = 0, \quad\quad \vartheta = z\frac{d}{dz}. 
\end{align*}
Then, as in (2.22)-(2.25) in \cite{KuMo}, we construct the functions
\begin{align*}
\phi_{1}(z) &= i e^{2\pi i\alpha} 
G_{3,0}^{0,3}\left(\left. \begin{array}{c} -\\ 0,-\alpha, - \alpha - \frac{1}{2}\end{array} \right| z e^{2\pi i}\right), \\
\phi_{2}(z) &= -i e^{-2\pi i\alpha} 
G_{3,0}^{0,3}\left(\left. \begin{array}{c} -\\ 0,-\alpha, - \alpha - \frac{1}{2}\end{array} \right| z e^{-2\pi i}\right), \\
\phi_{3}(z) &= 
G_{3,0}^{0,3}\left(\left. \begin{array}{c} -\\ 0,-\alpha, - \alpha - \frac{1}{2}\end{array} \right| z\right), \\
\phi_{4}(z) &= \phi_{1}(z)+\phi_{2}(z).
\end{align*}
Here the notation $z e^{2\pi i}$ means that we have analytically continued $\phi_3$ along a counterclockwise loop around the origin, and similarly for $z e^{-2\pi i}$. Then $\Psi$ is defined as
\begin{align} \label{eq:defPhiaaa}
\Psi(\zeta) = 
\begin{cases}
\begin{pmatrix} \phi_{1}(\zeta) & \phi_{2}(\zeta) & \phi_{3}(\zeta)\\ 
	\vartheta \phi_{1}(\zeta) & \vartheta \phi_{2}(\zeta) & \vartheta \phi_{3}(\zeta) \\
	\vartheta^2 \phi_{1}(\zeta) & \vartheta^2 \phi_{2}(\zeta) & \vartheta^2 \phi_{3}(\zeta)\end{pmatrix}, & 0 < \arg (\zeta)< \frac{\pi}{2}, \\
\begin{pmatrix} \phi_{4}(\zeta) & \phi_{2}(\zeta) & \phi_{3}(\zeta)\\ 
	\vartheta \phi_{4}(\zeta) & \vartheta \phi_{2}(\zeta) & \vartheta \phi_{3}(\zeta)\\
	\vartheta^2 \phi_{4}(\zeta) & \vartheta^2 \phi_{2}(\zeta) & \vartheta^2 \phi_{3}(\zeta)\end{pmatrix}, & \frac{\pi}{2} < \arg (\zeta)<\pi, \\
\begin{pmatrix}  \phi_{2}(\zeta) & -\phi_{1}(\zeta) & \phi_{3}(\zeta)\\ 
\vartheta \phi_{2}(\zeta) & -\vartheta \phi_{1}(\zeta) & \vartheta \phi_{3}(\zeta)\\  
\vartheta^2 \phi_{2}(\zeta) & -\vartheta^{2}\phi_{1}(\zeta) & \vartheta^2 \phi_{3}(\zeta)\end{pmatrix}, & -\frac{\pi}{2}<\arg(\zeta)<0,\\
\begin{pmatrix} \phi_{4}(\zeta) & -\phi_{1}(\zeta) & \phi_{3}(\zeta)\\ 
	\vartheta \phi_{4}(\zeta) & -\vartheta \phi_{1}(\zeta) & \vartheta \phi_{3}(\zeta)\\
	\vartheta^2 \phi_{4}(\zeta) & -\vartheta^2 \phi_{1}(\zeta) & \vartheta^2 \phi_{3}(\zeta)\end{pmatrix}, & -\pi<\arg (\zeta)<-\frac{\pi}{2}.
\end{cases}
\end{align}
The counterpart of $\Psi$ is $\Phi_\alpha$ as in (2.27) in \cite{KuMo}. In Lemma 3.4 in \cite{KuMo} it was shown that, as $\zeta\to\infty$
\begin{align} \label{psiAsymp}
\Psi(\zeta) = \left(\mathbb I + \frac{C_1}{\zeta} + \mathcal O\left(\frac{1}{\zeta^2}\right)\right) B(\zeta) e^{n \theta(\zeta)},
\end{align}
where $C_1$ is a $3\times 3$ matrix that depends only on $\alpha$, 
\begin{align*}
\theta(\zeta) &= -3 \zeta^\frac{1}{3} \times \begin{cases}
\begin{pmatrix}  \omega & 0 & 0\\ 0 & \omega^2  & 0\\ 0 & 0 & 1 \end{pmatrix}, & \operatorname{Im}(\zeta)>0, \\
\begin{pmatrix} \omega^2 & 0 & 0\\ 0 & \omega & 0\\ 0 & 0 & 1 \end{pmatrix}, & \operatorname{Im}(\zeta)<0,
\end{cases}\\
B(\zeta) &= \frac{2\pi}{\sqrt 3} T_\alpha^{-1} \zeta^{-\frac{2\beta}{3}}
\begin{pmatrix} \zeta^{-\frac{1}{3}} & 0 & 0 \\
		0 & 1 & 0 \\ 0 & 0 & \zeta^{\frac{1}{3}} \end{pmatrix}
		\times
	\begin{cases} \begin{pmatrix} 	\omega^2 & \omega & 1\\ 
	1 & 1 & 1\\ \omega & \omega^2 & 1 \end{pmatrix}  
	\begin{pmatrix} e^{\frac{2\pi i \beta}{3}} & 0 & 0 \\
	0 & e^{-\frac{2\pi i \beta}{3}} & 0 \\
	0 & 0 & 1 \end{pmatrix}, &\operatorname{Im}(z)>0,\\
	\begin{pmatrix}  \omega & -\omega^2 & 1\\ 1 & -1 & 1\\ 
	\omega^2 & -\omega & 1\end{pmatrix}
	\begin{pmatrix} e^{-\frac{2\pi i \beta}{3}} & 0 & 0 \\
	0 & e^{\frac{2\pi i \beta}{3}} & 0 \\
	0 & 0 & 1 \end{pmatrix},
	 &\operatorname{Im}(z)<0,
\end{cases}
\end{align*}
where $T_\alpha$ is a constant lower-triangular matrix that depends only on $\alpha$ (see (3.9) in \cite{KuMo}). The first factor on the right-hand side of \eqref{psiAsymp} can be extended to a full asymptotic series in negative powers of $\zeta$, in principle. $\Psi$ solves the following bare parametrix problem with constant jumps (see p. 17 in \cite{KuMo}). 
\begin{itemize}
\item[RH-$\Psi$1] $\Psi : \mathbb{C}\setminus (\mathbb R\cup i\mathbb R) \to \mathbb{C}^{3\times 3}$ is analytic. 
\item[RH-$\Psi$2] On $\mathbb R\cup i\mathbb R$ we have the following jumps (all curves oriented outwards).
\begin{align*} 
\Psi_{+}(\zeta) &= \Psi_{-}(\zeta)\begin{pmatrix} 0 & 1 & 0\\ -1 & 0 & 0\\ 0 & 0 & 1\end{pmatrix}, && \zeta\in \mathbb R^+, \\
\Psi_{+}(\zeta) &= \Psi_{-}(\zeta) \begin{pmatrix} 1 & 0 & 0\\ 1 & 1 & 0\\ 0 & 0 & 1\end{pmatrix}, && \zeta\in i\mathbb R^\pm,\\ 
\Psi_{+}(\zeta) &= \Psi_{-}(\zeta) \begin{pmatrix} 1 & 0 & 0\\ 0 & 0 & i e^{2\pi i\alpha}\\ 0 & -i e^{2\pi i\alpha} & 0\end{pmatrix}, && \zeta\in \mathbb R^-. 
\end{align*}
\item[RH-$\Psi$3] As $\zeta\to\infty$
\begin{align*}
\Psi(\zeta) = \left(\mathbb I + \frac{C_1}{\zeta} + \mathcal O\left(\frac{1}{\zeta^2}\right)\right) B(\zeta) e^{n \theta(\zeta)}.
\end{align*}
\item[RH-$\Psi$4] $\Psi D$ has the same behavior as $S$ has as $z\to 0$. 
\end{itemize}
As it turns out, a convenient choice for $\mathring E_n$ is given by
\begin{align*}
\mathring E_n(z) = n^{-2\beta} \left(\frac{f(z)}{z}\right)^{-\frac{2\beta}{3}} T_\alpha^{-1}.
\end{align*}
We could have avoided having to define $\mathring E_n$, with the $n^{-2\beta}$ factor in particular, if we included a factor $\zeta^\frac{2\beta}{3}$ in the definition of $\Psi(\zeta)$ and a factor $z^{-\frac{2\beta}{3}}$ in the definition of $D(z)$. We did not make this choice because we do not want to deviate too much from the expressions in \cite{KuMo}. Some straightforward algebra shows that  
\begin{align*}
\mathring P(z) N(z)^{-1} = \left(\mathbb I + \frac{T_\alpha^{-1} C_1 T_\alpha}{n^3 f(z)} + \mathcal O\left(\frac{1}{n^6 z^2}\right)\right) E(z)^{-1}
\end{align*}
uniformly on any circle that shrinks slower than order $n^{-3}$ as $n\to\infty$, where 
\begin{align*}
E(z) =n^{2\beta} \left(\frac{f(z)}{z}\right)^{\frac{2\beta}{3}} N(z) D(z)^{-1} e^{-n D_\varphi(z)-\theta(n^3 f(z))} B(n^3 f(z))^{-1} T_\alpha.
\end{align*}
It so happens (see (5.36) in \cite{KuMo}), that the $\varphi$-functions have an expansion in terms of powers of $z^\frac{1}{3}$ close to $z=0$. In particular, one can show (see Proposition 5.8 in \cite{KuMo}) that
\begin{align*}
n D_\varphi(z) + \theta(n^3 f(z)) = \mathcal O\left(n z^\frac{2}{3}\right)
\end{align*}
as $z\to 0$. This hints that we should try to obtain the matching on a shrinking circle of radius $n^{-\frac{3}{2}}$. A circle of fixed radius would lead to an exponential error. For $z\in \partial D(0,n^{-\frac{3}{2}})$ we have uniformly that
\begin{align*}
\mathring P(z) N(z)^{-1} = \left(\mathbb I + \frac{C(z)}{n^3 z} + \mathcal O\left(n^{-3}\right)\right) E(z)^{-1}
\end{align*}
as $n\to\infty$, where 
\begin{align*}
C(z) = \frac{z}{f(z)} T_\alpha^{-1} C_1 T_\alpha. 
\end{align*}
Indeed, $C$ is uniformly bounded on $\partial D(0,n^{-\frac{3}{2}})$ as $n\to\infty$ (in fact on $\overline{D(0,r)}$ entirely) and it is analytic. The function $E$ coincides with $E_n$ in (5.44) of \cite{KuMo}. Then there are some properties for $E_n$ that are proved in \cite{KuMo} (see Lemma 5.10(c) and Lemma 5.13 for the details), that we may immediately apply. Namely, $E$ is non-singular and analytic. But also, we have uniformly for $z\in\partial D(0,n^{-\frac{3}{2}})$ as $n\to\infty$ that
\begin{align*}
E(z) = \mathcal O(n)\quad\quad\text{and}\quad\quad E(z)^{-1} = \mathcal O(n),
\end{align*}
and uniformly for $z,w\in\partial D(0,n^{-\frac{3}{2}})$ as $n\to\infty$ that
\begin{align*}
E(z)^{-1} E(w) = \mathbb I + \mathcal O\left(n^\frac{5}{2}(z-w)\right).
\end{align*}
I suspect that the proof of Lemma 5.13 in \cite{KuMo} provides a general strategy to obtain the relevant estimates for $E$. Now we have all the ingredients to apply Theorem \ref{lem:matching}, setting $a=\frac{3}{2}, b=3, c=3, d=2, e=\frac{5}{2}$, provided that the inequalities $a\leq e<b$ and $d<\min(b,c)$ are valid, which is indeed the case. Hence we obtain a double matching where the jumps on both the inner and the outer circle behave as $\mathbb I + \mathcal O\left(\frac{1}{n}\right)$. One is now able to construct the final transformation to $R$, as in \eqref{finalT} (with an additional local parametrix around $q$). In this particular case we have
$$
\frac{a+c-e}{b-e} = 4,
$$
thus we have $K=1$ in Definition \ref{def:prefactors}. This means that $K+1=2$ iteration steps are needed in the construction of $E_n^0$. This is an improvement with respect to \cite{KuMo}, where (effectively) three iteration steps were needed to get the same degree of approximation. This discrepancy is caused by the fact that iteration steps in \cite{KuMo} got rid of only the term with the lowest order pole in each step. Then the amount of iteration steps needed grows linearly rather than logarithmically. 

We may apply Lemma \ref{thm:scalingb} to conclude that for $z$ of order $n^{-3}$ we have as $n\to\infty$
\begin{align*}
E_n^0(z) = E(0) \left(\mathbb I + \mathcal O(n^{-\frac{1}{2}})\right)
\quad\quad\text{and}\quad\quad
E_n^0(z)^{-1} = \left(\mathbb I + \mathcal O(n^{-\frac{1}{2}})\right) E(0)^{-1},
\end{align*}
and for $z,w$ of order $n^{-3}$ we have as $n\to\infty$
\begin{align*}
E_n^0(z)^{-1} E_n^0(w) = \mathbb I + \mathcal O(n^\frac{5}{2}(z-w)). 
\end{align*}
In \cite{KuMo} these estimates were used to obtain a scaling limit for an associated correlation kernel. Following the derivation in \cite{KuMo} (see Lemma 6.4), and using the notation
\begin{align*}
x_n = \frac{x}{f'(0) n^3}\quad\text{and}\quad y_n = \frac{y}{f'(0) n^3},
\end{align*}
we have for $n$ big enough that
\begin{multline*} 
	\frac{1}{n^3} K_n(x_n,y_n) = e^{\frac{2}{3}n (V(x_n) - V(y_n))} \\
	\frac{1}{2\pi i(x-y)} \begin{pmatrix} -1 & 1 & 0 \end{pmatrix}
	\Psi_+(n^3 f(y_n))^{-1}
	E_n^0(y_n)^{-1} R^{-1}(y) R(x) E_n^0(x_n)
	\Psi_+(n^3 f(x_n)) \begin{pmatrix} 1 \\ 1 \\ 0 \end{pmatrix}.
\end{multline*}
One can show with the help of Appendix A in \cite{BlKu} that we have the estimate \eqref{eq:estimateRtoI} with $\delta=1$. Remarkably, $\delta$ actually \textit{is} bigger than $a-b+d=\frac{1}{2}$, hence we may directly apply \eqref{eq:naive} to conclude that
\begin{align*} 
	\frac{1}{n^3} K_n(x_n,y_n) = 
	\frac{1}{2\pi i(x-y)} \begin{pmatrix} -1 & 1 & 0 \end{pmatrix}
	\Psi_+(y)^{-1} \Psi_+(x) 
	\begin{pmatrix} 1 \\ 1 \\ 0 \end{pmatrix} + \mathcal O\left(n^{-\frac{1}{2}}\right).
\end{align*}
uniformly for $x, y$ in compact sets as $n\to\infty$. 

Notice that, if we nevertheless wanted to apply Theorem \ref{thm:corKerScaling}, then we would need to take an extra term $C_2\zeta^{-2}$ in the expansion \eqref{psiAsymp}, i.e., we would then have $c=\frac{9}{2} \geq \frac{15}{4} = \frac{3}{2}a+2d-e$ and  
\begin{align*}
C(z) = \frac{z}{f(z)} \left(T_\alpha^{-1} C_1 T_\alpha + \frac{T_\alpha^{-1} C_2 T_\alpha}{n^3 z} \frac{z}{f(z)}\right).
\end{align*}
In an article that is under construction, we will show that the case $\theta=\frac{1}{r}$ with $r$ a positive integer, concerning an $(r+1)\times (r+1)$ RHP, corresponds to $a=\frac{r+1}{2}, b=r+1, d=r$ and $e = r+\frac{1}{2}$. In the case $r>2$, we are not so lucky that we can directly apply \eqref{eq:naive} and we really have to use Theorem \ref{thm:corKerScaling} instead. We omit the details as these will appear in the upcoming article.

\subsection{Cauchy-Laguerre three-chain}

In the $4\times 4$ RH analysis for the Cauchy-Laguerre three-chain in \cite{BeBo} Bertola and Bothner arrive at the following situation. In \cite{BeBo} one has a global parametrix $M$ (see Proposition 4.8 in \cite{BeBo}) and a local parametrix $Q$ (see (4.35) in \cite{BeBo}), which is of the form
\begin{align} \label{eq:defBeBoQ}
Q(z) = B_0(z) G^{(3)}(\zeta(z))
\left(\frac{2n}{3^\frac{3}{4}} e_1(z)\right)^{-A}
\left\{\begin{array}{ll}
e^{4\Omega \zeta^\frac{1}{4}(z) + \frac{n}{2} 3^{-\frac{3}{4}} e_2(z) \overline \Omega}, & \operatorname{Im}(z)>0,\\
e^{4\overline \Omega \zeta^\frac{1}{4}(z) + \frac{n}{2} 3^{-\frac{3}{4}} e_2(z) \Omega}, & \operatorname{Im}(z)>0.
\end{array}\right.
\end{align}
The factor on the far right, with $\Omega$ being a constant diagonal matrix, $e_1, e_2$ being some analytic functions and $\zeta(z) \approx \frac{16}{27} n^4 z$ being a conformal map around $0$ (see p. 1110 in \cite{BeBo} for these four expressions) , should be identified with the $e^{n D_\varphi}$ factor in \eqref{eq:defmathringP}. The factor in the middle, with $A$ some constant diagonal matrix (see (4.25) in \cite{BeBo}), should be identified with the $D$ factor in \eqref{eq:defmathringP}. The function $G^{(3)}$ (see (4.33) in \cite{BeBo}), called the bare Meijer-G parametrix, solves a corresponding bare parametrix problem and it should be identified with $\Psi$ in \eqref{eq:defmathringP}. It is constructed with the help of Meijer G-functions that solve a fourth order linear differential equation. Notice that $\zeta(z)$ has the interpretation of $n^b f(z)$ in \eqref{eq:defmathringP} with $b=4$ (indeed $n^{-4}\zeta$ is independent of $n$). The function $B_0$ (see (4.36) in \cite{BeBo}) is the first attempt for the analytic prefactor.

In \cite{BeBo} it is shown that $Q$ and $M$ satisfy a relation (see (4.38) in \cite{BeBo}) of the following asymptotic form. As $z\to\infty$ we have
\begin{align} \label{BeBoMatch}
Q(z) M(z)^{-1} \sim \widehat{M}(z) z^{-\frac{1}{8}\lambda_4} 
\left(\mathbb I+\sum_{j=1}^\infty K_j \zeta^{-\frac{j}{4}}\right)
H(z) z^{\frac{1}{8}\lambda_4} \widehat{M}(z)^{-1},
\end{align}
for $z\in \overline{D(0,r)}$ for a fixed $r>0$ sufficiently small. Here $\widehat{M}$ (see (4.27) in \cite{BeBo}) is $n$-independent and analytic in some neighborhood of $z=0$, $\lambda_4=\operatorname{diag}[3,1,-1,-3]$ (see (4.28) in \cite{BeBo}), $K_j$ are constant $n$-independent matrices and $H$ (see (4.39) in \cite{BeBo}) is some function with properties that we turn to in a second.

In \cite{BeBo} the matching was obtained with an iterative technique, in which one kept modifying $B_0$, which seems to be applicable to more general situations and larger sizes in principle. However, I suspect that the technical nature of the approach would make this a cumbersome procedure when the size of the RHP is large. For our purposes we can make the following identifications.
\begin{align} \nonumber
\mathring P(z) &= \left(\frac{z}{\zeta}\right)^{\frac{1}{8}\lambda_4} \widehat{M}(z)^{-1} Q(z),\\ \nonumber
N(z) &= M(z),\\ \nonumber
E(z) &= \widehat{M}(z) z^{-\frac{1}{8}\lambda_4} H(z)^{-1}  \zeta^{\frac{1}{8}\lambda_4}.
\end{align}
Comparing with \eqref{eq:defBeBoQ} and \eqref{eq:defmathringP} this then corresponds to the first guess of the prefactor in \eqref{eq:defmathringP} being
\begin{align*}
\mathring E_n(z) = \left(\frac{z}{\zeta}\right)^{\frac{1}{8}\lambda_4} \widehat{M}(z)^{-1} B_0(z).
\end{align*}
With these definitions we obtain by \eqref{BeBoMatch} the asymptotic expansion
\begin{align}\label{BeBoMatch2}
\mathring P(z) N(z)^{-1} E(z) \sim \zeta^{-\frac{1}{8}\lambda_4} 
\left(\mathbb I+\sum_{j=1}^\infty K_j \zeta^{-\frac{j}{4}}\right)
\zeta^{\frac{1}{8}\lambda_4}
\end{align}
as $\zeta\to 0$. In fact, a straightforward adaptation of the proof of Proposition 4.15 in \cite{BeBo} shows that the right-hand side of \eqref{BeBoMatch2} represents an asymptotic series in integer powers of $\zeta$ and this yields the situation \eqref{eq:matchingAsympS}. In particular, taking the asymptotic series up to $j=11$, we find ($n$-independent) constants $\tilde K_1$ and $\tilde K_2$ such that
\begin{align} \label{eq:BeBointegerAsymp}
\mathring P(z) N(z)^{-1} E(z) = \mathbb I + \frac{\tilde K_1}{\zeta} + \frac{\tilde K_2}{\zeta^2} + \mathcal O(\zeta^{-\frac{9}{4}}) 
= \mathbb I + \frac{\tilde K_1}{\zeta} + \mathcal O(\zeta^{-2})
\end{align}
as $\zeta\to 0$. The function $E$ is indeed analytic, as follows from Proposition 4.13 in \cite{BeBo}. For $z\in \partial D(0,n^{-\frac{4}{3}})$ it turns out, using (4.39) in \cite{BeBo}, that $H$ and $H^{-1}$ are uniformly bounded. Combining the boundedness of $H^{-1}$ and $H$ with the specific form of $\lambda_4$
  yields that uniformly for $z\in \partial D(0,n^{-\frac{4}{3}})$
\begin{align} \nonumber
E(z) = \mathcal O(n^\frac{3}{2}) \quad \text{ and }\quad E(z)^{-1} = \mathcal O(n^\frac{3}{2})
\end{align} 
as $n\to\infty$. Also, the boundedness of $H$ and $H^{-1}$, and the specific form of $\lambda_4$ yield that uniformly for $z_1,z_2\in \partial D(0,n^{-\frac{4}{3}})$ as $n\to\infty$
\begin{align} \nonumber
E(z_1)^{-1} E(z_2) &= \zeta(z_1)^{-\frac{1}{8}\lambda_4} H(z_1) 
z_1^{\frac{1}{8} \lambda_4} (\mathbb I + \mathcal O(z_1-z_2)) z_2^{-\frac{1}{8} \lambda_4}
H(z_2)^{-1} \zeta(z_2)^{\frac{1}{8} \lambda_4}\\ \nonumber
&= \zeta(z_1)^{-\frac{1}{8}\lambda_4} H(z_1) 
(\mathbb I + \mathcal O(n^\frac{4}{3} (z_1-z_2))+\mathcal O(n (z_1-z_2))) 
H(z_2)^{-1} \zeta(z_2)^{\frac{1}{8} \lambda_4}\\ \nonumber
&= \zeta(z_1)^{-\frac{1}{8}\lambda_4} 
(\mathbb I + \mathcal O(n^\frac{4}{3} (z_1-z_2)))
\zeta(z_2)^{\frac{1}{8} \lambda_4}\\ \nonumber
&= \mathbb I + \mathcal O(n^{\frac{4}{3}}(z_1-z_2)) + \mathcal O(n^{1+\frac{4}{3}+1}(z_1-z_2))\\ \nonumber
&=  \mathbb I + \mathcal O(n^{\frac{10}{3}}(z_1-z_2)).
\end{align}
In light of Theorem \ref{lem:matching} we identify $a=\frac{4}{3}, b=4, d=3$ and $e=\frac{10}{3}$. We still need to find a convenient value for $c$. If we take two terms of integer powers of $\zeta$ in the asymptotic expansion on the right-hand side of \eqref{BeBoMatch2} then this seems to be enough, i.e., we identify
\begin{align} \nonumber
C(z) = \frac{n^4 z}{\zeta(z)} \tilde K_1.
\end{align}
Notice that this will indeed make $C$ uniformly bounded on $\partial D(0,n^{-\frac{4}{3}})$ as $n\to\infty$ (and in fact on any circle inside the domain of $\zeta$). Indeed, using \eqref{eq:BeBointegerAsymp} we now have uniformly for $z\in\partial D(0,n^{-\frac{4}{3}})$ that
\begin{align}\label{BeBoMatch3}
\mathring P(z) N(z)^{-1} E(z) = \mathbb I + \frac{C(z)}{n^4 z} + \mathcal O(\zeta^{-2})
=\mathbb I+ \frac{C(z)}{n^4 z}+\mathcal O(n^{-\frac{16}{3}})
\end{align}
as $n\to\infty$, and we identify $c=\frac{16}{3}$. 
Then Theorem \ref{lem:matching} provides us with $E_n^0, E_n^\infty$ such that as $n\to\infty$
\begin{align} \nonumber
E_n^0(z) \mathring P(z) &= \left(\mathbb I + \mathcal O(n^{-\frac{7}{3}})\right) E_n^\infty(z) N(z), &\text{uniformly for }z\in \partial D(0,n^{-a}),\\ \nonumber
E_n^\infty(z) &= \mathbb I + \mathcal O(n^{-1}), &\text{uniformly for }z\in \partial D(0,r).
\end{align}
In order to use Theorem \ref{thm:corKerScaling} we should have that $c\geq \min(\frac{3}{2}a+d,\frac{3}{2} a+2d-e) = \frac{14}{3}$. This is indeed the case, so we do not need to consider more terms in the asymptotic expansion \eqref{BeBoMatch2}. 

It should be relatively easy to extend the result in \cite{BeBo} to $p$-chains with $p>3$, the construction of a corresponding global parametrix for general $p$ being the only missing ingredient. 


\pagebreak

\subsection{Non-intersecting squared Bessel paths: critical case}

In \cite{KuMFWi} an issue with the matching was also encountered in a $3\times 3$ RHP. There, Kuijlaars, Mart\'inez-Finkelshtein and Wielonsky resolved the problem by modifying the global parametrix. This example is a bit special in that the corresponding functions in \eqref{eq:defmathringP} are not entirely what one would expect, as we shall see. It is a nice insight that Theorem \ref{lem:matching} can nevertheless be applied here. The global parametrix $N_\alpha$ can be found in (4.9) in \cite{KuMFWi}. 
The local parametrix in \cite{KuMFWi} (see (5.13)) has the form
\begin{align*}
Q(z) = E_n(z) \Phi_\alpha\left(n^\frac{3}{2} f(z); n^\frac{1}{2} \tau(z)\right) 
\operatorname{diag}[1,z^\alpha,1]
\operatorname{diag}\left[e^{-n\lambda_1(z)}, e^{-n \lambda_2(z)}, e^{-n \lambda_3(z)}\right]
e^{\frac{2 n z}{3 t(1-t)}}.
\end{align*}

The diagonal factor on the right (see (2.22) in \cite{KuMFWi} for the $\lambda$ functions) should be identified with $e^{n D_\varphi}$ from \eqref{eq:defmathringP} and the diagonal factor in the middle should be identified with $D$ from \eqref{eq:defmathringP}. Perhaps it is best to absorb $e^{\frac{2 n z}{3 t(1-t)}}$ into the definition of $e^{n D_\varphi}$, we will make this choice. $E_n$ (see (5.44) in \cite{KuMFWi}) could be interpreted as the first guess for the analytic prefactor in \eqref{eq:defmathringP}, we will define $\mathring E$ differently though. The function $f$ is a conformal map and $\tau$ is some analytic function (see (5.40) and Lemma 5.5 in \cite{KuMFWi}). In \cite{KuMFWi} one has a bare parametrix problem with a solution $\Phi_\alpha(z; \tau)$ (see Definition 5.1 in \cite{KuMFWi}) for some parameter $\alpha$, which we identify with $\Psi$ in \eqref{eq:defmathringP}. It is constructed with the help of solutions to the linear differential equation $z \phi'''+\alpha \phi'' - \tau \phi'-\phi=0$.  One finds the corresponding asymptotic behavior using a classical steepest descent analysis. In particular, in \cite{KuMFWi} (see Lemma 5.3) it is shown that as $z\to\infty$
\begin{align} \nonumber
\Phi_\alpha(z; \tau) 
= \frac{i}{\sqrt 3} L_\alpha(z) \times
\left\{\begin{array}{ll}
\left(\mathbb I + \frac{M_\alpha^+(\tau)}{z^\frac{1}{3}}+\mathcal O(z^{-\frac{1}{3}})\right)
\begin{pmatrix}
e^{\theta_1(z,\tau)} & 0 & 0\\
0 & e^{\theta_2(z,\tau)} & 0\\
0 & 0 & e^{\theta_3(z,\tau)}
\end{pmatrix} & \operatorname{Im}(z)>0,\\
\left(\mathbb I + \frac{M_\alpha^-(\tau)}{z^\frac{1}{3}}+\mathcal O(z^{-\frac{2}{3}})\right)
\begin{pmatrix}
e^{\theta_2(z,\tau)} & 0 & 0\\
0 & e^{\theta_1(z,\tau)} & 0\\
0 & 0 & e^{\theta_3(z,\tau)}
\end{pmatrix} & \operatorname{Im}(z)<0,
\end{array}\right.
\end{align}
where $\theta_k(z,\tau)= \frac{3}{2} \omega^{2k} z^\frac{2}{3}+\tau \omega^k z^\frac{1}{3}$ (see (1.21) in \cite{KuMFWi}), $\omega=e^\frac{2\pi i}{3}$, $M^\pm_\alpha(\tau)$ (see (5.26) and (5.27) in \cite{KuMFWi}) depends polynomially on $\tau$ and (see (5.25) in \cite{KuMFWi})
\begin{align*}
L_\alpha(z) = z^{-\frac{\alpha}{3}} z^{\frac{1}{3}\operatorname{diag}[1,0,-1]} \times
\left\{\begin{array}{rl}
\begin{pmatrix} \omega & \omega^2 & 1\\ 1 & 1 & 1\\ \omega^2 & \omega & 1\end{pmatrix} \operatorname{diag}[e^{\alpha\pi i/3}, e^{-\alpha\pi i/3},1], & \operatorname{Im}(z)>0,\\
\begin{pmatrix} \omega^2 & -\omega & 1\\ 1 & -1 & 1\\ \omega & -\omega^2 & 1\end{pmatrix} \operatorname{diag}[e^{-\alpha\pi i/3}, e^{\alpha\pi i/3},1], & \operatorname{Im}(z)<0.
\end{array}\right.
\end{align*}
In fact the expression $\mathbb I + M_\alpha^\pm(\tau) z^{-\frac{1}{3}}+\mathcal O(z^{-\frac{2}{3}})$ comes from an asymptotic series
\begin{align*}
\mathbb I + \frac{M_{\alpha,1}^\pm(\tau)}{z^\frac{1}{3}} + \frac{M_{\alpha,2}^\pm(\tau)}{z^\frac{2}{3}} + \ldots
\end{align*}
in $z^{-\frac{1}{3}}$, for certain coefficients $M_{\alpha,1}^\pm(\tau), M_{\alpha,2}^\pm(\tau),\ldots$ (these are not to be found in \cite{KuMFWi}). Of course we should then identify $M_{\alpha,1}^\pm(\tau)=M_\alpha^\pm(\tau)$. The authors note that $\Phi_\alpha$ solves a specific first-order matrix-valued ODE. From this ODE one can deduce that $(M_{\alpha,k}^\pm(\tau))_k$ satisfies a three term recurrence relation $M^\pm_{\alpha,k+2}(\tau)=\mathcal A^\pm M_{\alpha,k+1}^\pm(\tau)+\mathcal B^\pm M^\pm_{\alpha,k}(\tau)$, with coefficients $\mathcal A^\pm, \mathcal B^\pm$ that are at most linear in $\tau$. Then we infer that there exists a constant $\mathfrak c>0$ such that  $\lvert\lvert M^\pm_{\alpha,k}(\tau) \rvert\rvert \leq \mathfrak c^k (1+ \tau^k)$ for all $\tau\in\mathbb C$. In particular, the $M_{\alpha,k}^\pm(\tau)$ are bounded when $\tau$ is bounded. 

The behavior of $\theta_1,\theta_2,\theta_3$ is a little extraordinary, in that they are not multiples of $z^\frac{3}{2}$, and that they are in fact functions of two variables $z$ and $\tau$. Nevertheless, Theorem \ref{lem:matching} can be applied here. Let us define the auxiliary fuction
\begin{align*}
F_\alpha(z,\tau) := \frac{i}{\sqrt 3} L_\alpha(z)
 \times
\left\{\begin{array}{ll}
\begin{pmatrix}
e^{\theta_1(z,\tau)} & 0 & 0\\
0 & e^{\theta_2(z,\tau)} & 0\\
0 & 0 & e^{\theta_3(z,\tau)}
\end{pmatrix} & \operatorname{Im}(z)>0,\\
\begin{pmatrix}
e^{\theta_2(z,\tau)} & 0 & 0\\
0 & e^{\theta_1(z,\tau)} & 0\\
0 & 0 & e^{\theta_3(z,\tau)}
\end{pmatrix} & \operatorname{Im}(z)<0.
\end{array}\right.
\end{align*}
As one can check (use for this the proof of Lemma 5.7 in \cite{KuMFWi}), $F_\alpha$ has the same jumps as $\Phi_\alpha$ has on the positive and negative ray (although not the jumps on the lenses). This means that $\Phi_\alpha(z; \tau) F_\alpha(z,\tau)^{-1}$ does not have jumps for $z$ on the real line. Plugging in the asymptotic expansion for $\Phi_\alpha$ we get
\begin{align} \label{eq:PhiPsiAsymp}
\Phi_\alpha(z; \tau) F_\alpha(z,\tau)^{-1} 
\sim L_\alpha(z) \left(\mathbb I + \frac{M_{\alpha,1}^\pm(\tau)}{z^\frac{1}{3}} + \frac{M_{\alpha,2}^\pm(\tau)}{z^\frac{2}{3}} + \ldots\right) L_\alpha(z)^{-1},
\end{align}
as $z\to\infty$. Since there is no jump on the real line, the expression on the right-hand side of \eqref{eq:PhiPsiAsymp} actually equals an asymptotic series in powers of $z^{-1}$. The coefficients of that asymptotic series are then also bounded when $\tau$ is bounded. We identify, with $E_n, Q, N_\alpha, \lambda_k, f(z), \tau(z)$ as in \cite{KuMFWi},
\begin{align*}
\mathring P(z) &= -\sqrt{3} i n^{\frac{\alpha}{2}} E_n(z)^{-1} Q(z),\\
N(z) &= N_\alpha(z),\\
E(z) &= \frac{i}{\sqrt 3} n^{-\frac{\alpha}{2}} N_\alpha(z) e^{-\frac{3 n z}{3 t(1-t)}}
\begin{pmatrix}
e^{n\lambda_1(z)} & 0 & 0\\
0 & z^{-\alpha} e^{n\lambda_2(z)} & 0\\
0 & 0 & e^{n \lambda_3}
\end{pmatrix}
F_\alpha(n^\frac{3}{2} f(z),n^\frac{1}{2}\tau(z))^{-1} \\
&= -\sqrt 3 i n^{-\frac{\alpha}{2}}  N_\alpha(z) \operatorname{diag}[1,z^{-\alpha},1] L_\alpha(n^\frac{3}{2} f(z))^{-1} = n^{-\frac{\alpha}{2}} E_n(z).
\end{align*}
Effectively, this amounts to putting $\mathring E_n$ equal to $-\sqrt 3 i n^\frac{\alpha}{2}$ in \eqref{eq:defmathringP}. We used the formula under (6.16) in \cite{KuMFWi} to rewrite $E$, that is,
\begin{align*}
\theta_k(n^\frac{3}{2} f(z); n^\frac{1}{2}\tau(z)) &= n \lambda_k(z) - \frac{2nz}{3 t(1-t)}, & k=1,2,3. 
\end{align*}
In this special case the authors actually managed to match the $e^{n D_\varphi}$ factor with the $e^\theta$ factor, as in \eqref{eq:behavPsi}, exactly. The cost is that they needed two functions $f$ and $\tau$ to make it work. It follows from Lemma 5.5 in \cite{KuMFWi} that $n^\frac{1}{2} \tau(z)$ is uniformly bounded on $\partial D(0,n^{-\frac{1}{2}})$. Then some straightforward algebra, using \eqref{eq:PhiPsiAsymp}, yields that
\begin{align} \label{eq:integerAsymp2}
\mathring P(z) N(z)^{-1} E(z) \sim L_\alpha(n^\frac{3}{2} f(z)) 
\left(\mathbb I + \frac{M_{\alpha,1}^\pm(n^\frac{1}{2}\tau(z))}{(n^\frac{3}{2} f(z))^\frac{1}{3}} + \frac{M_{\alpha,2}^\pm(n^\frac{1}{2}\tau(z))}{(n^\frac{3}{2} f(z))^\frac{2}{3}} + \ldots\right)
 L_\alpha(n^\frac{3}{2} f(z))^{-1},
\end{align}
uniformly for $z\in\partial D(0,n^\frac{1}{2})$ as $n\to\infty$, where the coefficients $M_{\alpha,k}^\pm(n^\frac{1}{2} \tau(z))$ are uniformly bounded on $\partial D(0,n^{-\frac{1}{2}})$. As we pointed out earlier, this is actually an asymptotic series with integer powers evaluated in $(n^\frac{3}{2} f(z))^{-1}$. We read off that we should take $a=\frac{1}{2}$ and $b=\frac{3}{2}$. It is not hard to see that
\begin{align} \label{eq:KuMFWiLa}
L_\alpha(n^\frac{3}{2} f(z)) = n^{-\frac{\alpha}{2}} \operatorname{diag}[n^\frac{1}{2},1,n^{-\frac{1}{2}}] L_\alpha(f(z)).
\end{align}
A simple adaptation of Lemma 5.7 in \cite{KuMFWi} yields that the function
\begin{align*}
 -\sqrt 3 i  N_\alpha(z) \operatorname{diag}[1,z^{-\alpha},1] L_\alpha(f(z))^{-1}
\end{align*}
defines an analytic function that does not depend on $n$. Combining this with \eqref{eq:KuMFWiLa} yields that
\begin{align*}
E(z) = \mathcal O(n^\frac{1}{2}) \quad\quad \text{ and }\quad \quad E(z)^{-1} = \mathcal O(n^\frac{1}{2})
\end{align*}
as $n\to\infty$, uniformly for $z\in\partial D(0,n^{-\frac{1}{2}})$. Thus we take $d= 1$. Similarly, we have uniformly for $z,w\in\partial D(0,n^{-\frac{1}{2}})$ that
\begin{align*}
E(z)^{-1} E(w) = \mathbb I + \mathcal O(n (z-w))
\end{align*}
as $n\to\infty$. Hence we take $e = 1$. It remains to find out how many terms we should take in the asymptotic expansion, i.e., what $c$ should be. In order to apply Theorem \ref{thm:corKerScaling} we need that ${c\geq \min(\frac{3}{2} a+d,\frac{3}{2} a+2d-e) = \frac{7}{4}}$. Then we should take
\begin{align*}
C(z) = \frac{z}{f(z)} \tilde M_1(n^\frac{1}{2}\tau(z)),
\end{align*}
where $\tilde M_1(n^\frac{1}{2}\tau(z))/ (n^\frac{3}{2} f(z))$ is defined as the second term in the asymptotic expansion in integer powers of $(n^\frac{3}{2} f(z))^{-1}$ in \eqref{eq:integerAsymp2}. The corresponding $c$ then equals $2$. As we mentioned before, $n^\frac{1}{2} \tau(z)$ is uniformly bounded on $\partial D(0,n^{-\frac{1}{2}})$ as $n\to\infty$, thus the coefficients $M_{\alpha,k}^\pm(n^\frac{1}{2} \tau(z))$, and consequently $\tilde M_1(n^\frac{1}{2}\tau(z))$, are also bounded. Hence $C$ is indeed uniformly bounded on $\partial D(0,n^{-\frac{1}{2}})$. We conclude that the assumptions of Theorem \ref{lem:matching} are satisfied. Indeed, by the maximum modulus principle we also have that $C$ is uniformly bounded on $\partial D(0,n^{-\frac{2}{3}})$, as we should have in order to apply Theorem \ref{thm:corKerScaling}. 

\vspace{1cm}

\Large \noindent\textbf{Acknowledgements} \normalsize \vspace{0.25cm}

\noindent The author is supported by a PhD fellowship of the Flemish Science Foundation (FWO) and partly by
the long term structural funding-Methusalem grant of the Flemish Government. The author is grateful to Arno Kuijlaars for proofreading the article and suggested improvements. The author would like to thank Thomas Bothner, Guilherme Silva and Alexander Tovbis for useful discussions.

\phantomsection

\end{document}